\documentclass[11pt]{amsart}
\usepackage{etex}
\usepackage[hmargin=3cm,vmargin=3cm]{geometry}

\usepackage[latin1]{inputenc}
\usepackage{xspace,amssymb,amsfonts, mathtools}
\usepackage{amsthm,amsmath,amscd,stmaryrd,latexsym}
\usepackage{palatino, euscript, newclude}
\usepackage[normalem]{ulem}
\usepackage[all,cmtip]{xy}

\usepackage{arydshln,pinlabel,hyperref}

\usepackage[nohug]{diagrams}
\usepackage[T1]{fontenc}
\makeatletter

\usepackage{color}
\usepackage{tikz, tikz-cd}
\usetikzlibrary{matrix,arrows,decorations.pathmorphing, decorations.markings}

\tikzstyle directed=[postaction={decorate,decoration={markings,
    mark=at position #1 with {\arrow{>}}}}]
\tikzstyle rdirected=[postaction={decorate,decoration={markings,
    mark=at position #1 with {\arrow{<}}}}]

\usepackage{comment}

\newtheorem{theorem}[equation]{Theorem}

\newtheorem{lemma}[equation]{Lemma}
\newtheorem{proposition}[equation]{Proposition}
\newtheorem{corollary}[equation]{Corollary}

\theoremstyle{definition}
\newtheorem{definition}[equation]{Definition}

\theoremstyle{remark}

\newtheorem{remark}[equation]{Remark}

\numberwithin{equation}{section}

\DeclareMathAlphabet{\mathbbold}{U}{bbold}{m}{n}

\newcommand{\ZZ}{\mathbb{Z}}
\newcommand{\NN}{\mathbb{N}}
\newcommand{\QQ}{\mathbb{Q}}
\newcommand{\CC}{\mathbb{C}}
\newcommand{\Bn}{\mathfrak{B}_n}
\newcommand{\HH}{\text{HH}}
\newcommand{\B}{\mathcal{B}}
\newcommand{\BB}{\tilde{\B}}
\newcommand{\T}{\mathcal{T}}
\newcommand{\FT}{\mathcal{FT}}
\renewcommand{\b}{\beta}
\newcommand{\holim}{\text{holim}\,}
\newcommand{\cone}{\mathrm{Cone}}

\newcommand{\mm}{\mathbf{m}}

\newcommand{\webC}{\mathbf{NWeb}}
\newcommand{\foamC}{\mathbf{NFoam}}
\newcommand{\qbinom}[2]{{#1 \brack #2}}

\newcommand{\sln}{\mathfrak{sl}_N}
\newcommand{\slinf}{\mathfrak{sl}_\infty}
\newcommand{\qsln}{\mathcal{U}_q(\sln)}
\newcommand{\C}[1]{C_N\left( #1 \right)}
\newcommand{\F}{\mathcal{F}}
\renewcommand{\H}{\mathcal{H}}

\newcommand{\AS}{\mathbf{A}}
\newcommand{\BS}{\mathbf{B}}
\newcommand{\CS}{\mathbf{C}}

\newcommand{\arc}{\,\tikz[baseline=\the\dimexpr\fontdimen22\textfont2\relax ]{\draw[->] plot [smooth] coordinates {(0,0)(0.05,0.2)(0,0.4)}; }\, }

\newcommand{\posxT}{\,\tikz[baseline=3mm]{
\draw[->] (0.6,0)--(0,1);
\draw[ultra thick, white] (0,0)--(0.6,1);
\draw[ultra thick, white] (0.05,0)--(0.65,1);
\draw[ultra thick, white] (-0.05,0)--(0.55,1);
\draw[->] (0,0)--(0.6,1);
}\,}

\newcommand{\posx}{\,\tikz[baseline=\the\dimexpr\fontdimen22\textfont2\relax]{
\draw[->] (0.4,0)--(0,0.4);
\draw[ultra thick, white] (0,0)--(0.4,0.4);
\draw[ultra thick, white] (0.05,0)--(0.45,0.4);
\draw[ultra thick, white] (-0.05,0)--(0.35,0.4);
\draw[->] (0,0)--(0.4,0.4);
}\,}

\newcommand{\posxij}{\,\tikz[baseline=3mm]{
\draw[->] (0.6,0)--(0,1);
\draw[ultra thick, white] (0,0)--(0.6,1);
\draw[ultra thick, white] (0.05,0)--(0.65,1);
\draw[ultra thick, white] (-0.05,0)--(0.55,1);
\draw[->] (0,0)--(.6,1);
\node at (0,-.2) {$_i$};
\node at (0.6,-.2) {$_j$};
\node at (0,1.2) {$_j$};
\node at (0.6,1.2) {$_i$};
}\,}

\newcommand{\negx}{\,\tikz[baseline=\the\dimexpr\fontdimen22\textfont2\relax]{
\draw[->] (0,0)--(0.4,0.4);
\draw[ultra thick, white] (0.4,0)--(0,0.4);
\draw[ultra thick, white] (0.4,0.05)--(0.05,0.4);
\draw[ultra thick, white] (0.4,0-0.05)--(-0.05,0.4);
\draw[->] (0.4,0)--(0,0.4);
}\,}

\newcommand{\negxT}{\,\tikz[baseline=3mm]{
\draw[->] (0,0)--(.6,1);
\draw[ultra thick, white] (.6,0)--(0,1);
\draw[ultra thick, white] (0.65,0)--(0.05,1);
\draw[ultra thick, white] (.55,0)--(-.05,1);
\draw[->] (0.6,0)--(0,1);
}\,}

\newcommand{\negxij}{\,\tikz[baseline=\the\dimexpr\fontdimen22\textfont2\relax]{
\draw[->] (0,0)--(0.6,1);
\draw[ultra thick, white] (0.6,0)--(0,1);
\draw[ultra thick, white] (0.6,0.05)--(0.05,1);
\draw[ultra thick, white] (0.6,0-0.05)--(-0.05,1);
\draw[->] (0.6,0)--(0,1);
\node at (0,-.2) {$_i$};
\node at (0.6,-.2) {$_j$};
\node at (0,1.2) {$_j$};
\node at (0.6,1.2) {$_i$};
}\,}

\newcommand{\ladder}[5]{\,\tikz[baseline=3mm]{
\draw[->] (0,0)--(0.0,1);
\draw[->] (0.8,0)--(0.8,1);
\draw[->] (0,.2)--(.8,.35);
\draw[->] (.8,.65)--(0,.8);
\node at (0,-.2) {$_{#1}$};
\node at (0.8,-.2) {$_{#2}$};
\node at (0,1.2) {$_{#3}$};
\node at (0.8,1.2) {$_{#4}$};
\node at (.4,.05) {$_{#5}$};
}\,}

\newcommand{\posxnegx}{\,\tikz[baseline=3mm]{
\draw (0,0)--(0.4,0.4);
\draw[ultra thick, white] (0.4,0)--(0,0.4);
\draw[ultra thick, white] (0.4,0.05)--(0.05,0.4);
\draw[ultra thick, white] (0.4,-0.05)--(-0.05,0.4);
\draw (0.4,0)--(0,0.4);
\draw[->] (0.4,0.4)--(0,0.8);
\draw[ultra thick, white] (0,0.4)--(0.4,0.8);
\draw[ultra thick, white] (0.05,0.4)--(0.45,0.8);
\draw[ultra thick, white] (-0.05,0.4)--(0.35,0.8);
\draw[->] (0,0.4)--(0.4,0.8);
}\,}

\newcommand{\arcs}{\,\tikz[baseline=\the\dimexpr\fontdimen22\textfont2\relax ]{\draw[->] plot [smooth] coordinates {(0,0)(0.05,0.2)(0,0.4)};\draw[->] plot [smooth] coordinates {(0.4,0)(0.35,0.2)(0.4,0.4)}; }\, }

\newcommand{\arcsT}{\,\tikz[baseline=3mm]{\draw[->] plot [smooth] coordinates {(0,0)(0.05,0.5)(0,1)};\draw[->] plot [smooth] coordinates {(0.4,0)(0.35,0.5)(0.4,1)}; }\, }

\newcommand{\ladderC}{\,\tikz[baseline=3mm ]{
\draw[->] (0,0)--(0.0,1);
\draw[->] (0.5,0)--(0.5,1);
\draw[fill=black]  (-.05,.4) rectangle (.55,.6);
 }\, }
 
\newcommand{\sladderC}{\,\tikz[baseline=\the\dimexpr\fontdimen22\textfont2\relax ]{
\draw[->] (0,0)--(0.0,0.4);
\draw[->] (0.2,0)--(0.2,0.4);
\draw[fill=black]  (-.025,.15) rectangle (.225,.25);
 }\, }
 
\newcommand{\posxcl}{\,\tikz[baseline=\the\dimexpr\fontdimen22\textfont2\relax]{
\draw[->] (0.4,0)--(0,0.4);
\draw[ultra thick, white] (0,0)--(0.4,0.4);
\draw[ultra thick, white] (0.05,0)--(0.45,0.4);
\draw[ultra thick, white] (-0.05,0)--(0.35,0.4);
\draw (0,0)--(0.4,0.4);
\draw plot [smooth] coordinates {(0.4,0.4)(0.5,0.2)(0.4,0)};
}\,}

\newcommand{\negxcl}{\,\tikz[baseline=\the\dimexpr\fontdimen22\textfont2\relax]{
\draw (0,0)--(0.4,0.4);
\draw[ultra thick, white] (0.4,0)--(0,0.4);
\draw[ultra thick, white] (0.4,0.05)--(0.05,0.4);
\draw[ultra thick, white] (0.4,0-0.05)--(-0.05,0.4);
\draw[->] (0.4,0)--(0,0.4);
\draw plot [smooth] coordinates {(0.4,0.4)(0.5,0.2)(0.4,0)};
}\,}

\newcommand{\barbell}[5]{\,\tikz[baseline=3mm]{
\draw[->-] (0,0)--(0.2,0.3);
\draw[->-] (0.4,0)--(0.2,0.3);
\draw[->-] (0.2,.3)--(.2,.7);
\draw[->] (.2,.7)--(0,1);
\draw[->] (.2,.7)--(.4,1);
\node at (0,-.2) {$_{#1}$};
\node at (0.4,-.2) {$_{#2}$};
\node at (0,1.2) {$_{#3}$};
\node at (0.4,1.2) {$_{#4}$};
\node at (.55,.5) {$_{#5}$};
}\,}

\newcommand{\mergeWeb}
{\,\tikz[baseline=\the\dimexpr\fontdimen22\textfont2\relax]{
\draw (0,0)--(.2,.2);
\draw (.4,0)--(.2,.2);
\draw (.2,.2)--(.2,.4);
}\,}

\newcommand{\splitWeb}
{\,\tikz[baseline=\the\dimexpr\fontdimen22\textfont2\relax]{
\draw (.2,0)--(.2,.2);
\draw (.2,.2)--(0,.4);
\draw (.2,.2)--(.4,.4);
}\,}

\tikzset{->-/.style={decoration={
  markings,
  mark=at position .5 with {\arrow{>}}},postaction={decorate}}}
  
\tikzset{-<-/.style={decoration={
  markings,
  mark=at position .5 with {\arrow{<}}},postaction={decorate}}}

\newcommand{\xymat}[1]{\entrymodifiers={+!!<0pt,\fontdimen22\textfont2>}\xymatrix{#1}}

\makeatletter
\DeclareRobustCommand{\rvdots}{%
  \vbox{
    \baselineskip4\p@\lineskiplimit\z@
    \kern-\p@
    \hbox{.}\hbox{.}\hbox{.}
  }}
\makeatother

\begin{document}

\title{colored Khovanov-Rozansky homology for infinite braids}
\author[Michael Abel and Michael Willis]{Michael Abel and Michael Willis}
\address{Department of Mathematics, Duke University, Durham, NC 27708}
\email{maabel \char 64 math.duke.edu}
\address{Department of Mathematics, UCLA, 520 Portola Plaza, Los Angeles, CA 90095}
\email{mike.willis \char 64 math.ucla.edu}

\begin{abstract}
We show that the limiting unicolored $\mathfrak{sl}(N)$ Khovanov-Rozansky chain complex of any infinite positive braid categorifies a highest-weight projector. This result extends an earlier result of Cautis categorifying highest-weight projectors using the limiting complex of infinite torus braids. Additionally, we show that the results hold in the case of colored HOMFLY-PT Khovanov-Rozansky homology as well. An application of this result is given in finding a partial isomorphism between the HOMFLY-PT homology of any braid positive link and the stable HOMFLY-PT homology of the infinite torus knot as computed by Hogancamp.
\end{abstract}

\maketitle

\section{Introduction}

The Jones-Wenzl projector is a special idempotent element of the Temperley-Lieb algebra representing a highest-weight projector in the representation theory of $\mathcal{U}_q(\mathfrak{sl}_2)$.  One approach to categorifying this projector involves the Khovanov complex of the infinite full twist $C_2(\FT_n^\infty)$, first studied by Rozansky \cite{Roz10a} who proved it was homotopy equivalent to the categorified Jones-Wenzl projector of Cooper-Krushkal \cite{CK12a}.  (For alternative, inequivalent approaches, see \cite{Kh05} and \cite{FSS} and later \cite{Rob} for the case of $\mathcal{U}_q(\mathfrak{sl}_3)$.)  Then in \cite{IW16}, Islambouli and the second author proved that the stable Khovanov complex of any positive semi-infinite complete braid $\BB$ on $n$ strands also converges to Cooper-Krushkal's categorification of the Jones-Wenzl projector via comparison with Rozansky's infinite twist complex.

Later work used the above idea of Rozansky in \cite{Roz10a} to study categorified highest weight projectors in other settings. Rose in \cite{Rose12} showed that the infinite full twist categorifies the highest-weight projector for tensor products of fundamental representations of $\mathcal{U}_q(\mathfrak{sl}_3)$. Also using this idea, Cautis in \cite{Cau12} proved the same result for infinite twists categorifying highest weight projectors for representations of $\mathcal{U}_q(\mathfrak{sl}_N)$ for all $N$.  
These highest weight projectors are the analogues of the Jones-Wenzl projector for $\mathcal{U}_q(\mathfrak{sl}_N)$; their categorification plays a central role in Cautis' extension of Khovanov-Rozansky homologies to non-fundamental representations, and may also be used to study the 2-representations of $\mathcal{U}_q(\mathfrak{sl}_N)$ in a similar manner to how the projectors themselves can be used to study the usual representations of $\mathcal{U}_q(\mathfrak{sl}_N)$.  All of this was also extended to the case of Khovanov-Rozansky HOMFLY-PT homology by Hogancamp \cite{Hog15} via computing the Rouquier complex of the infinite full twist, a semi-infinite chain complex of bimodules in the categorification of the Hecke algebra.  Other limits of the infinite full twist in Khovanov and Khovanov-Rozansky homology can be used to categorify more general projectors, as in the work of Rozansky \cite{Roz10b} and of the first author with Hogancamp \cite{AbHog15a}.

In this text, we explore the generalization of the work in \cite{IW16} to the case of Khovanov-Rozansky $\mathfrak{sl}(N)$ and HOMFLY-PT homology. We may color our braid by elements of $\{1,...,N-1\}$ for the $\mathfrak{sl}(N)$ case, and by positive integers in $\NN$ for the HOMFLY-PT case respectively. We denote a braid $\b$ in which all strands are colored by the same integer $m$ by $\b_{(m)}$. Our main result is the following theorem.

\begin{theorem}[See Theorems \ref{thm-infBraid} and \ref{thm-infBraidHHH}]\label{thm-infBraidintro}
	Let $\BB$ be a complete semi-infinite positive braid on $n$ strands.  Then there exists a well-defined Khovanov-Rozansky $\mathfrak{sl}(N)$ complex $\C{\BB_{(m)}}$ for all $m \in \{1,...,N-1\}$.  Furthermore, 
	$$\C{\BB_{(m)}} \simeq \C{\FT_{(m)}^\infty}.$$
	Likewise, the Khovanov-Rozansky HOMFLY-PT complex $\F(\BB_{(m)})$ exists for all $m \in \NN$ and
	$$\F(\BB_{(m)}) \simeq \F(\FT_{(m)}^\infty).$$
	
	
\end{theorem}

We will clarify the notion of a complete semi-infinite positive braid, along with the definitions of the Khovanov-Rozansky complexes, in \S \ref{subsec-foams} and \S \ref{subsec-HHH}. The proof of this theorem is similar to the proof of the main theorem of \cite{IW16}. A rough outline is as follows. Let $\B_\ell$ denote the finite braid consisting of the first $\ell$ crossings of a semi-infinite braid word $\B$ representing $\BB$. $\C{\B_{\ell,(m)}}$ can be written as a mapping cone containing $\C{\FT_{(m)}^{z(\ell)}}$ as a quotient complex, for some positive integer $z(\ell) > 0$. This number $z(\ell)$ roughly counts the number of full twists present in $\B_\ell$. We prove that the natural projection map $F_{\ell}: \C{\B_{\ell,(m)}} \to \C{\FT_{(m)}^{z(\ell)}}$ is an isomorphism up to a certain homological degree $|F_{\ell}|_h$. We prove that $|F_{\ell}|_h \to \infty$ as $\ell \to \infty$ by carefully keeping track of shifts in homological degree by pulling and twisting `forks' and `ladders' through crossings. Unlike in the case of pulling turnbacks in Khovanov homology, these moves do not always lower the number of crossings in the diagram and so we will need some more careful combinatorics.  The proof in the HOMFLY-PT case is almost identical, with only a few minor differences.

Although the statement of Theorem \ref{thm-infBraidintro} is for positive semi-infinite complete braids, the results can be extended to more general infinite braids which are bi-infinite, containing a finite number of negative crossings, or non-complete. In the first two cases, we obtain the same complex up to a grading shift. In the latter case, we obtain a diagram with a finite number of crossings and multiple infinite full twists side-by-side. We refer the reader to \S \ref{subsec-mainResult} for more precise statements. 

One immediate application of the main theorem in the HOMFLY-PT case is given in studying the homology of braid positive links. Recent work of Elias, Hogancamp, and Mellit \cite{EH16,Hog17,Mel17} gives an explicit computation of the Khovanov-Rozansky HOMFLY-PT homology of torus links. Using our bounds on the homological order $|F_{\ell}|_h$ we are able to give a partial isomorphism between the HOMFLY-PT homology of a braid positive link and the stable HOMFLY-PT homology of torus links. In the following theorem $\H$ is used to denote Khovanov-Rozansky HOMFLY-PT homology.

\begin{theorem}[See Theorem \ref{thm-linkEst}]  \label{thm-linkEstIntro}
Let $\b$ be a finite positive braid on $n$ strands with $\ell$ crossings and suppose $y(\ell)$ is the number of ``diagonals'' in $\b$ (see Definition \ref{def-diagsAndZones}) . Then there exists a map $\bar{F}_\ell: \H(\b) \to \H(\T_n^{y(\ell)})$ which is an isomorphism for all homological degrees less than $y(\ell)$. 	
\end{theorem}

\subsection{Outline of paper} In Section 2 we will review the relevant background needed for colored Khovanov-Rozansky $\mathfrak{sl}(N)$ complexes of infinite braids, as well as recalling the results of Cautis from \cite{Cau12}. In Section 3 we restate our main results and corollaries and give a detailed proof. In Section 4 we explore the case of colored Khovanov-Rozansky HOMFLY-PT homology and discuss the estimation theorem for braid positive links.

\subsection*{Acknowledgements}
The authors would like to thank Sabin Cautis, David Rose and Paul Wedrich for helpful conversations and Mikhail Khovanov for asking a question which eventually led to this work.  They would also like to thank an anonymous referee whose comments and suggestions very much improved the paper.  The first author was supported in part by NSF grant DMS-1406371; the second author was supported in part by NSF grant DMS-1612159.

\section{Background}
In this section we will recall relevant background material from the Queffelec-Rose foam construction of colored Khovanov-Rozansky homology. We will then describe the diagrammatic lemmas and tools from homological algebra we will use to prove our main result. Finally we will recall a special case of our main result, proved by Cautis, in terms of infinite full twists. Though our result subsumes this one in certain cases, it is an important key piece in proving our result.

\subsection{The $\mathfrak{sl}(N)$ foam category and Khovanov-Rozansky homology}\label{subsec-foams} Colored Khovanov-Rozansky homology was first constructed independently by Wu \cite{Wu09} and Yonezawa \cite{Yo11}. This homology theory generalizes the original construction of Khovanov and Rozansky \cite{KR08} and categorifies the colored $SL(N)$ polynomial when coloring components by fundamental representations. Queffelec and Rose gave a combinatorial/geometric construction of colored Khovanov-Rozansky homology in terms of ``webs'' and ``foams'' \cite{QR16}. It is this construction which we will briefly recall here.

We begin with the category $\webC$, as described by Cautis, Kamnitzer, and Morrison in \cite{CKM}.  The objects of $\webC$ are given by sequences $\mm = (m_1,...,m_\ell)$ where $\ell > 0$ and $m_i \in \{0,1,...,N\}$ together with a zero object.  The 1-morphisms are formal sums of upward oriented trivalent graphs with edges labeled by integers in the same coloring set $\{0,1,...,N\}$.  At any vertex, the labels of the two edges of similar orientation (incoming or outgoing) are required to sum up to give the label of the third edge.  Such graphs are generated by the basic graphs in Figure \ref{fig-basicWebs}.  There is also a set of local relations for such trivalent graphs; we list a few of the most important ones in Figure \ref{fig-webRel}, and refer the reader to \cite{CKM} for a complete list.

Note that $\webC$ has a monoidal structure where $\mm_1 \boxtimes \mm_2$ is given by concatenation of sequences and for two webs $\Gamma_1$ and $\Gamma_2$, $\Gamma_1\boxtimes \Gamma_2$ is given by horizontal juxtaposition. That is, if $\Gamma_1: \mm_1 \to \mm_1'$ and $\Gamma_2: \mm_2 \to \mm_2'$ then $\Gamma_1 \boxtimes \Gamma_2 : \mm_1 \boxtimes \mm_2 \to \mm_1' \boxtimes \mm_2'$.

\begin{figure}[ht]
\begin{tikzpicture}
	\draw[thick,->-] (0,0)--(0,1); 
	\node at (0,-.25) {$i$};
	\node at (0,1.25) {$i$};
	\draw[thick,->-] (2,0)--(2,.5); \node at (2,-.25) {$i+j$};
	\draw[thick,->-] (2,.5)--(1.5,1); \node at (1.5,1.25) {$i$};
	\draw[thick,->-] (2,.5)--(2.5,1);\node at (2.5,1.25) {$j$};
	\draw[thick,->-] (4,0)--(4.5,.5); \node at (4,-.25) {$i$};
	\draw[thick,->-] (5,0)--(4.5,.5); \node at (5,-.25) {$j$};
	\draw[thick,->-] (4.5,.5)--(4.5,1); \node at (4.5,1.25) {$i+j$};
\end{tikzpicture}
\caption{Basic $\sln$-webs \label{fig-basicWebs}}	
\end{figure}

\begin{figure}[ht]
\begin{tikzpicture}
\draw[thick,->-] (2,0)--(2,.5); \node at (2,-.25) {$i+j$};
	\draw[thick,->-] (2,.5)--(1.5,1); \node at (1.25,1.25) {$i$};
	\draw[thick,->-] (2,.5)--(2.5,1);\node at (2.75,1.25) {$j$};
	\draw[thick,->-] (1.5,1)--(2,1.5); 
	\draw[thick,->-] (2.5,1)--(2,1.5); 
	\draw[thick,->-] (2,1.5)--(2,2);\node at (2,2.25) {$i+j$};
	\node[scale=1.4] at (3.75,1) {$= \qbinom{i+j}{i}$};
	\draw[thick,->-] (5,0)--(5,2); \node at (5,-.25) {$i+j$};\node at (5,2.25) {$i+j$};
\end{tikzpicture}
\\
\begin{tikzpicture}
	\draw[thick,->-] (0,0)--(0,0.6);\draw[thick,] (0,0.6)--(0,1.4);\draw[thick,->] (0,1.4)--(0,2);
	\draw[thick,->-] (1,0)--(1,0.6);\draw[thick,] (1,0.6)--(1,1.4);\draw[thick,->] (1,1.4)--(1,2);
	\draw[thick,-<-] (0,.75)--(1,.45);\draw[thick,->-] (0,1.25)--(1,1.55);
	\node at (0,-.25) {$i$};\node at (1,-.25) {$j+k$};\node at (.5,.4) {$\ell$};\node at (0,2.25) {$i+k$};\node at (1,2.25) {$j$};
	\node[scale=1.4] at (3.65,1) {$= \sum_{p=\max(0,k)}^{\ell} \qbinom{j-i}{\ell-p}$};
	\draw[thick,->-] (6.25,0)--(6.25,0.6);\draw[thick,] (6.25,0.6)--(6.25,1.4);\draw[thick,->] (6.25,1.4)--(6.25,2);
	\draw[thick,->-] (7.25,0)--(7.25,0.6);\draw[thick,] (7.25,0.6)--(7.25,1.4);\draw[thick,->] (7.25,1.4)--(7.25,2);
	\draw[thick,->-] (6.25,.45)--(7.25,.75);\draw[thick,-<-] (6.25,1.55)--(7.25,1.25);
	\node at (6.25,-.25) {$i$};\node at (7.25,-.25) {$j+k$};\node at (6.75,1.7) {$p$};\node at (6.25,2.25) {$i+k$};\node at (7.25,2.25) {$j$};
\end{tikzpicture}
\[[n] = \dfrac{q^n-q^{-n}}{q-q^{-1}}, \quad [n]! = [n][n-1]\cdots[1],\quad \qbinom{n}{k} = \dfrac{[n]!}{[n-k]![k]!}\]
\caption{Selected relations for $\sln$-webs.  See \cite{CKM} for a full list. \label{fig-webRel}}
\end{figure}

We interpret any such graph as a mapping between the sequence of labels at the bottom to the sequence of labels at the top. These graphs are called \emph{$\sln$-webs} due to their relation to the representation theory of $\qsln$. Sometimes we will omit edges labeled by $0$ and $N$, but allowing these labels in the definition will make later definitions easier to write. By convention, we will also allow edges labeled by integers larger than $N$ for the sake of later definitions. However, any web with such an edge will be set equal to the zero web (that is, the unique morphism factoring through the zero object). 

We now move on to the $2$-category $\foamC$.  The objects and 1-morphisms of $\foamC$ are the same as those for $\webC$, except that the 1-morphisms are considered as direct sums rather than genuine sums, and the relations between such 1-morphisms are discarded for the moment.  The 2-morphisms are matrices of labeled singular cobordisms between $\sln$-webs, called $\sln$-foams. These cobordisms are generated by the basic cobordisms in Figure \ref{fig-basicFoams}.

\begin{figure}[ht]
\include*{BasicFoams}
\caption{The basic generating $\sln$-foams in $\foamC$ (Image from \cite{QR16}).  \label{fig-basicFoams}}
\end{figure}

Similar to the convention for $\sln$-webs, we will interpret $\sln$-foams as mapping from the bottom boundary to the top boundary. Each facet of an $\sln$-foam will be labeled with an element of $\{0,1,...,N\}$.  Any facet whose boundary is shared with an edge of a web must share the same label as that edge of the web.  We allow decorations $\bullet_p$ on the facets of the foams where $p$ is a symmetric polynomial in a number of variables equal to the label of the facet.  There also exists a set of local relations for these $\sln$-foams which allow for a lifting of the web relations in $\webC$ (such as those in Figure \ref{fig-webRel}) to 2-isomorphisms between the corresponding 1-morphisms in $\foamC$.  The reader should consult \S3 of \cite{QR16} for more details. 

The category $\foamC$ admits a grading. We will note the gradings in cases that it is necessary, but will once again refer the reader to \cite{QR16} for more complete information. All chain complexes of foams are assumed to have degree 0 differentials. We will denote grading shifts in $\foamC$ with the notation $q^k$ for a grading shift upwards by $k$. 

To any tangle diagram $T$ whose components are labeled by elements of $\{0,1,...,N\}$, we can associate a chain complex in $\foamC$, which we will denote by $C_N(T)$. The homotopy equivalence type of $C_N(T)$ is an isotopy invariant of the tangle $T$. In diagrams and figures, we will often omit the notation $C_N(\bullet)$ and simply draw the corresponding tangle or web unless there is a chance for confusion. In this text we will exclusively focus on the case that the tangle $T$ is actually a braid or braid closure.

To build a chain complex in $\foamC$ for a braid, we construct basic chain complexes for each crossing. Suppose that $i \leq j$, then

\begin{equation}\label{eq-posxDef}
\xymat{
\posxij :=& \ladder{i}{j}{j}{i}{0}\ar[r]^-{d_0} & tq\ladder{i}{j}{j}{i}{1} \ar[r]^-{d_1}& \cdots \ar[r]^-{d_{i-1}}& t^iq^i\ladder{i}{j}{j}{i}{i}
}
\end{equation}

\begin{equation}\label{eq-negxDef}
\xymat{
\negxij:=& \ladder{i}{j}{j}{i}{i}\ar[r]^-{d'_{i-1}}&\cdots \ar[r]^-{d'_{1}}& t^{i-1}q^{i-1}\ladder{i}{j}{j}{i}{1}
\ar[r]^-{d_{0}'}&t^iq^i\ladder{i}{j}{j}{i}{0}	}.
\end{equation}

We use the variable $t$ to denote a shift in homological degree with the convention that all webs naturally sit in homological degree 0. In this paper, we use the conventions that differentials have homological degree $1$. We remark that the labeled edges determine all edges in each web, and that certain webs may be zero webs if they have a label larger than $N$.

The maps $d_k$ are degree $0$ foams as specified in \cite{QR16}. The maps $d_k'$ are the same foams, but reflected to switch the source and target webs. Finally, if $i > j$, then we reflect each web around a vertical axis and perform the analogous transformation to the foams $d_k$ and $d'_k$.

We note that our conventions differ from those of \cite{QR16}. In particular, our $C_N(\posx)$ differs from their definition of $C_N(\posx)$ by a shift of $t^iq^i$ (where the two strands are labeled by $i$ and $j$ and $i \leq j$). However, our convention makes studying stabilization behavior more straightforward. We remark that under our convention,

\begin{equation}\label{eq-shiftRMoves}
 \posxnegx \simeq t^iq^i \arcs, \quad \posxcl \simeq q^{i(i-N)} \arc , \quad \negxcl \simeq t^iq^{i(N-i+1)} \arc \,\,\,\text{ (All edges labeled by $i$)} .
\end{equation}

That is, Reidemeister I and Reidemeister II moves hold only up to a shift. Since Reidemeister III moves only change the arrangement of crossings rather than their number (or orientation), it is not hard to see that applying Reidemeister III moves does not introduce any new shifts.

Finally, we record a few key observations which we will use in the sequel. Suppose that both strands of $\posx$ are labeled by the same integer $m \in \{1,...,N-1\}$. Then, we can write $\C{\posx}$ as 
\begin{equation}\label{eq-uniclrCplx}
	\begin{tikzcd}
	\posxT = \arcsT \arrow{r}{d_0} & tq \ladder{\,}{\,}{\,}{\,}{1} \arrow{r}{d_1} & \cdots \arrow{r}{d_{m-2}}& t^{m-1}q^{m-1} \ladder{\,}{\,}{\,}{\,}{m-1} \arrow{r}{d_{m-1}} &t^mq^m\ladder{\,}{\,}{\,}{\,}{m}.
	\end{tikzcd}
\end{equation}
In the above diagram, some of the webs may be the zero web depending on $N$ and $m$.

\begin{remark}\label{rmk-0edges}
In the final term in the complex (\ref{eq-uniclrCplx}) we have a ladder whose rungs are labeled by $m$, and with an edge labelled by $0$. By our earlier convention on deleting edges, this indicates a barbell shape $$\barbell{\,}{\,}{\,}{\,}{2m}:=\ladder{\,}{\,}{\,}{\,}{m}.$$  We will choose to include or omit the edge labeled by $0$ as convenient for later formulae and descriptions.  See Remark \ref{rmk-0edgesAgain}. 
\end{remark}
Define a subcomplex of $\posx$, which we will denote by $\sladderC$, by

\begin{equation}\label{eq-blkboxDef}
\begin{tikzcd}
\ladderC = tq\ladder{\,}{\,}{\,}{\,}{1} \arrow{r}{d_1}& \cdots \arrow{r}{d_{m-2}}& t^{m-1}q^{m-1} \ladder{\,}{\,}{\,}{\,}{m-1} \arrow{r}{d_{m-1}}& t^mq^m\ladder{\,}{\,}{\,}{\,}{m}.
\end{tikzcd}	
\end{equation}

We will often call the webs in \sladderC ``ladders''. The following lemma combines (\ref{eq-uniclrCplx}) and (\ref{eq-blkboxDef}). 

\begin{lemma}\label{lem-posxCone}
$C_N(\posx)$ is homotopy equivalent to a mapping cone 

\begin{equation}\label{eq-posxCone}
\posxT \simeq \left[\,\,\arcsT \xrightarrow{d_0} t^{-1} \ladderC \right]\,\,. 	
\end{equation}

Note, we use the nonstandard notation $[A \xrightarrow{f} B]$ to denote the mapping cone of $f$. Likewise, $C_N(\negx)$ is homotopy equivalent to a mapping cone 

\begin{equation}\label{eq-negxCone}
\negxT \simeq \left[\,\,\ladderC \xrightarrow{d'_0} t^{m-1}q^m \arcsT \right]\,\,. 	
\end{equation}
\end{lemma}

In particular, Lemma \ref{lem-posxCone} gives us a canonical projection
\begin{equation}\label{eq-posxProj}
\pi: \C{\posx} \longrightarrow \C{\arcs}
\end{equation}
whose mapping cone is homotopy equivalent to \sladderC.  Our constructions will be based upon this projection $\pi$.

\subsection{Fork sliding and twisting}\label{subsec-forkST} We now introduce some diagrammatic relations on complexes of $\sln$-webs which will be important in the sequel. These relations follow from similar relations in \cite{QR16}, \cite{Yo11} and \cite{Wu09}, but differ slightly due to our conventions.

\begin{remark}\label{rmk-0edgesAgain}
Note that the formulas presented below hold even if one of the strands is labelled by a 0, which indicates that the strand is in fact not present.  See Remark \ref{rmk-0edges} above.
\end{remark}

\begin{figure}[ht]
\includegraphics[scale=.35]{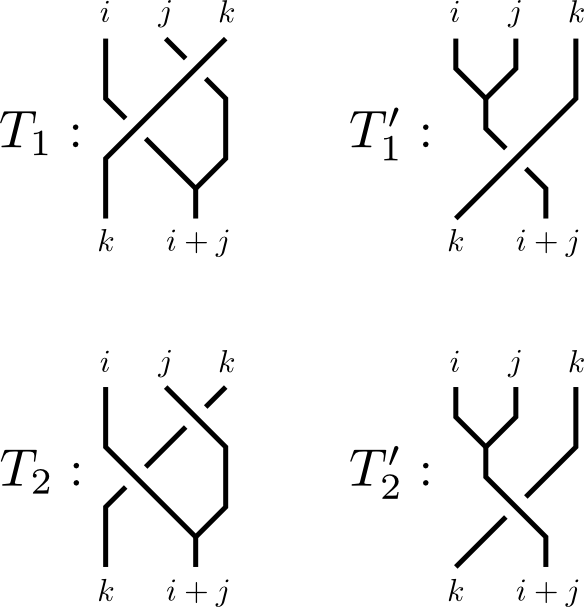}
\caption{The colored tangles $T_1,T_1',T_2,T_2'$ referenced in Lemma \ref{lem-forkSlide}.  The phrase "fork sliding" refers to sliding the vertex (the forking point) under or over another strand.}
\label{fig-forkSlide}
\end{figure}

\begin{lemma}[Fork Sliding] \label{lem-forkSlide}
	Let $T_1, T_1',T_2,T_2'$ be the colored tangles in Figure \ref{fig-forkSlide}.  Then 
	\begin{equation}\label{eq-forkSlide}
		\begin{aligned}
			C_N(T_1) &\simeq (tq)^{\min(i,k)+\min(j,k) - \min(i+j,k)}C_N(T_1')\\
			C_N(T_2) &\simeq C_N(T_2')
		\end{aligned}
	\end{equation}
\end{lemma}

Lemma \ref{lem-forkSlide} follows immediately from the invariance of fork sliding in the conventions of \cite{QR16}, after accounting for the shift in our definition of $C_N(\posx)$. We note that $\min(i,k)+\min(j,k) \geq \min(i+j,k)$ and that the inequality is strict when $i = j = k$.

\begin{figure}[ht]
\includegraphics[scale=.35]{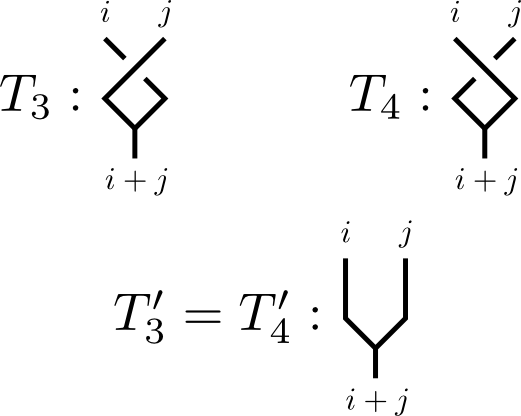}
\caption{The colored tangles $T_3,T_4$ and $T_3'=T_4'$ referenced in Lemma \ref{lem-forkTwist}.  The phrase "fork twisting" refers to a Reidemeister 1-like twist applied at the vertex (the forking point).}
\label{fig-forkTwist}
\end{figure}

\begin{lemma}[Fork Twisting] \label{lem-forkTwist}
	Let $T_3, T_3',T_4,T_4'$ be the colored tangles in Figure \ref{fig-forkTwist}.  Then 
	\begin{equation}\label{eq-forkTwist}
		\begin{aligned}
			C_N(T_3) &\simeq t^{\min(i,j)}q^{ij+\min(i,j)}C_N(T_3')\\
			C_N(T_4) &\simeq q^{-ij}C_N(T_4')
		\end{aligned}
	\end{equation}
\end{lemma}

Once again, the proof of Lemma \ref{lem-forkTwist} follows from the isomorphism for fork twisting in \cite{QR16} after taking into account our shift for $C_N(\posx)$. Often we will be performing similar operations on ladder webs as well. We record these isomorphisms in the following two lemmas where all crossings are positive. The proofs and statements for other cases are similar, but not used in this text.

\begin{figure}[ht]
\includegraphics[scale=.35]{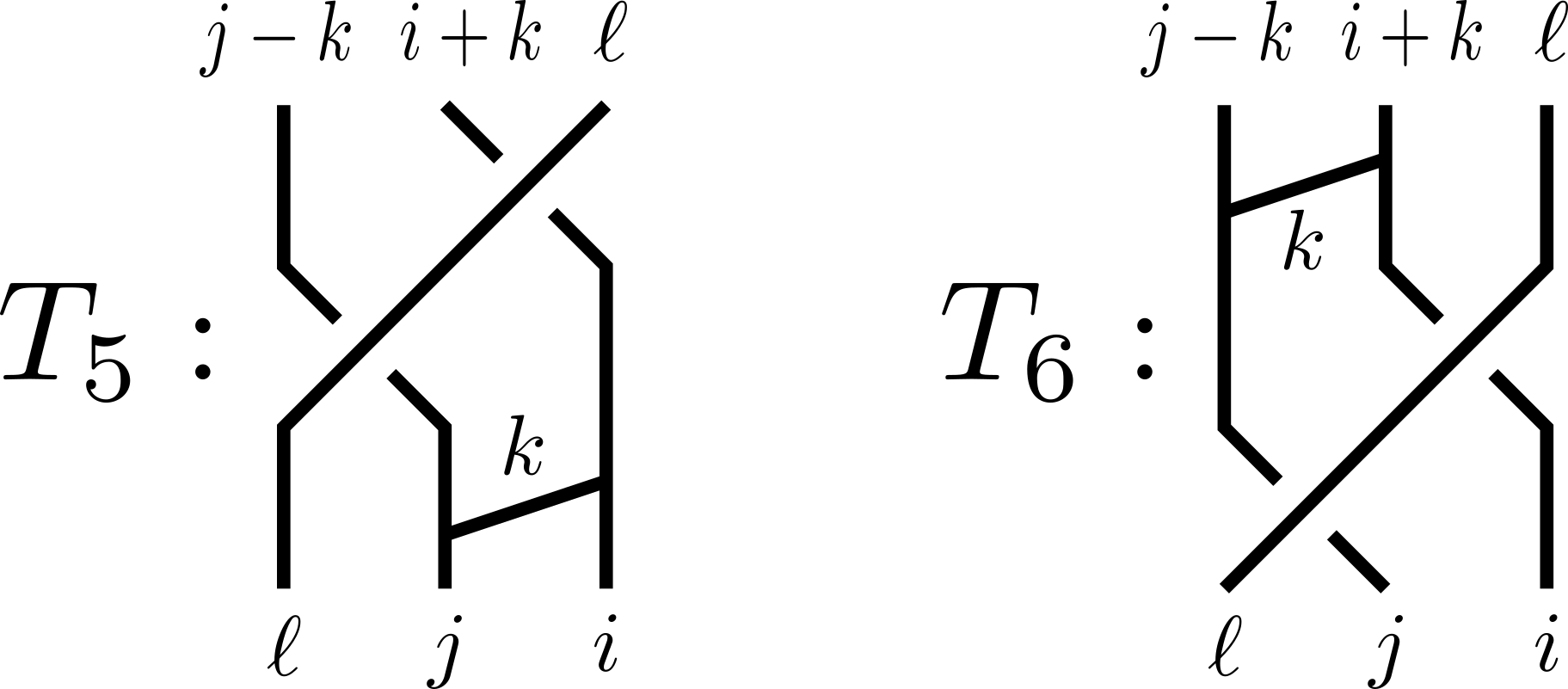}
\caption{The colored tangles $T_5$ and $T_6$ referenced in Lemma \ref{lem-ladSlide}.  The phrase "ladder sliding" refers to sliding the "rung" (here colored by $k$) under another strand.  There are similar results for sliding a rung over another strand.}
\label{fig-ladSlide}
\end{figure}

\begin{lemma}[Ladder sliding]\label{lem-ladSlide}
	Let $T_5$ and $T_6$ be the colored tangle diagrams in Figure \ref{fig-ladSlide}. Then $\C{T_5}\simeq(tq)^\alpha\C{T_6},$ where $\alpha = \min(i,\ell) + \min(j,\ell) - \min(i+k,\ell) - \min(j-k,\ell)$.
\end{lemma}

\begin{proof}
The isomorphism comes from a composition of two fork slides, illustrated in Figure \ref{fig-ladSlideProof}:
\begin{figure}[ht]
\includegraphics[scale=.35]{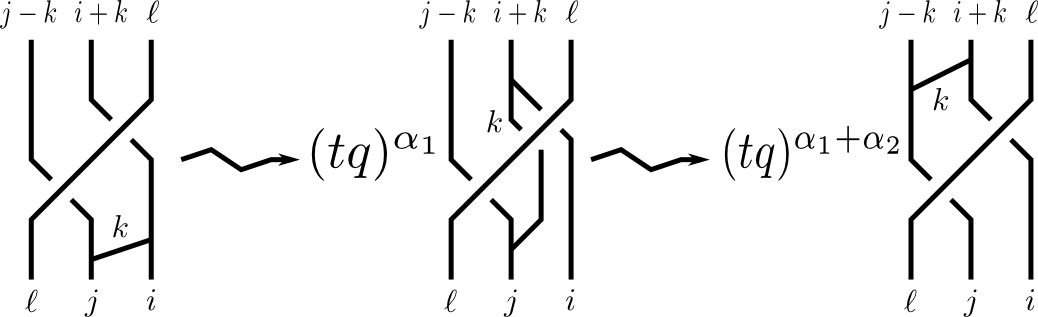}
\caption{The two fork slides that are utilized to complete a ladder slide.  The $t$ and $q$ indicate homological and $q$-degree shifts respectively.}
\label{fig-ladSlideProof}
\end{figure}
By Lemma \ref{lem-forkSlide}, $\alpha_1 = \min(i+k,\ell)-\min(i,\ell)-\min(k,\ell)$ and $\alpha_2 = \min(k,\ell)+\min(j-k,\ell) - \min(j,\ell)$. Thus, $\alpha = \alpha_1 + \alpha_2 = \min(i,\ell) + \min(j,\ell) - \min(i+k,\ell) - \min(j-k,\ell)$ as desired.
\end{proof}

\begin{figure}[ht]
\includegraphics[scale=.35]{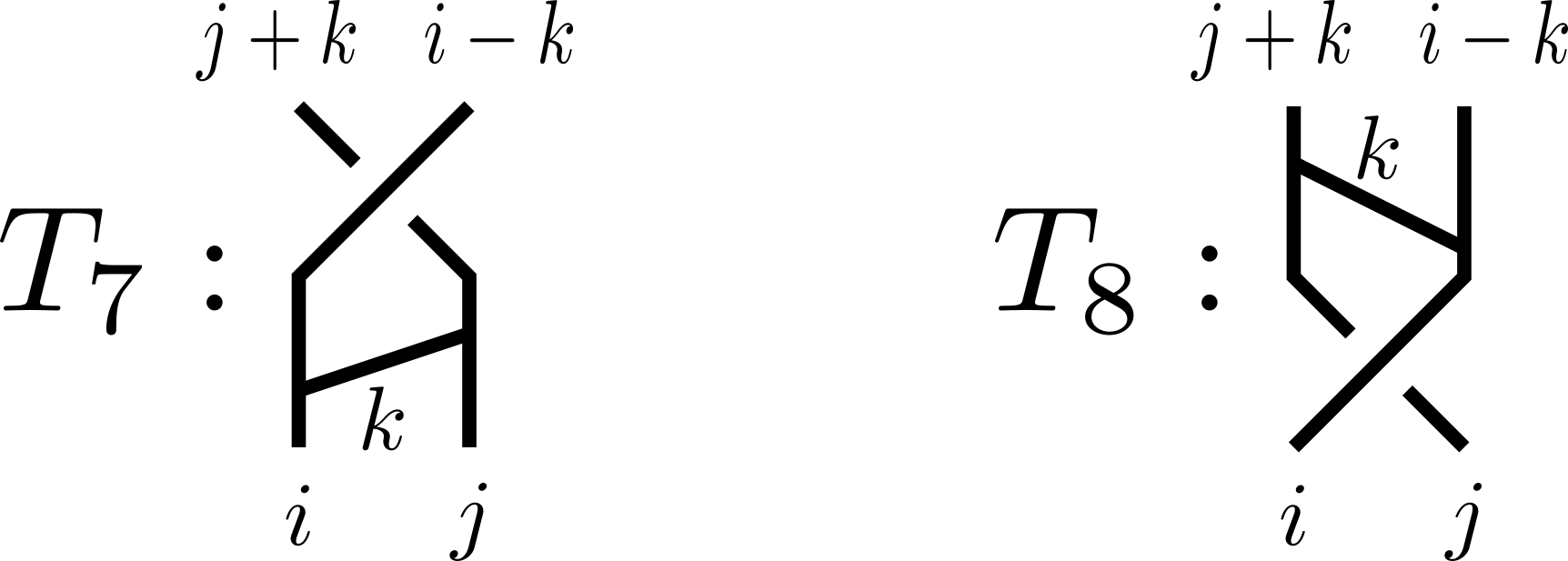}
\caption{The colored tangles $T_7$ and $T_8$ in Lemma \ref{lem-ladTwist}.  The phrase "ladder twisting" refers to the $180^\circ$ twist applied to the "rung" of the ladder as it is "pulled through" the crossing.}
\label{fig-ladTwist}
\end{figure}

\begin{lemma}[Ladder twisting]\label{lem-ladTwist}
Let $T_7$ and $T_8$ be the colored tangle diagrams in Figure \ref{fig-ladTwist}, then $$\C{T_7} \simeq t^{\min(i-k,j-k)-\min(i,k)}q^{\min(i-k,j-k)-\min(i,k)+(i-j-k)k}\C{T_8}.$$
\end{lemma}

\begin{proof} The isomorphism is realized by the composition of fork slides and fork twists illustrated in Figure \ref{fig-ladTwistProof}.
\begin{figure}[ht]
\includegraphics[scale=.35]{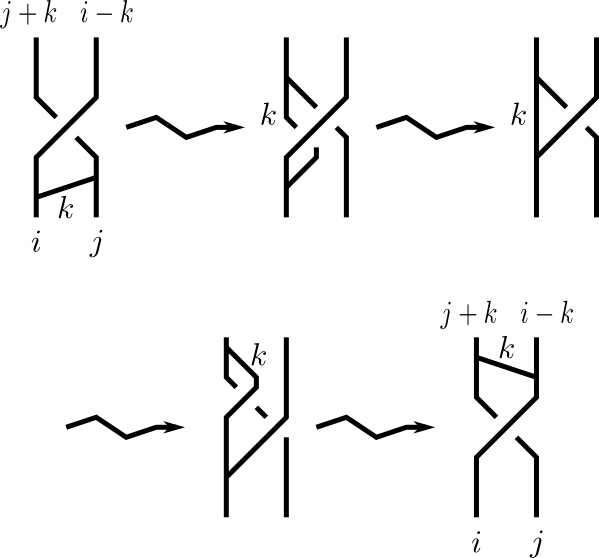}
\caption{The sequence of fork slides and fork twists that are utilized to complete a ladder twist.  For simplicity, the degree shifts are omitted, and during the intermediate stages of the move only the placement of the color $k$ is indicated.}
\label{fig-ladTwistProof}
\end{figure}
The first move is a fork slide, and so by Lemma \ref{lem-forkSlide} we see an initial shift of $$(qt)^{\min(i-k,j+k)-\min(k,i-k)-\min(j,i-k)}.$$
The second move untwists a fork, so that Lemma \ref{lem-forkTwist} gives the shift $$(qt)^{\min(i-k,k)}q^{(i-k)k}.$$
The third move introduces a fork twist, so Lemma \ref{lem-forkTwist} now gives a shift of $$(qt)^{-\min(j,k)}q^{-jk}.$$
Finally, we perform one more fork slide giving a fourth shift $$(qt)^{\min(j,k)+\min(j,i-k)-\min(i,j)}.$$ The total grading shift indicated in the statement of the lemma is the product of all four of these shifts. 
\end{proof}

All of the shifts indicated in the lemmas above may seem a bit unwieldy at first, but if we focus only on the homological shifts then they can be summarized tidily in a single proposition.  First, we introduce a bit of notation.

\begin{definition}\label{def-crossingMin}
For a given crossing $\chi$ between two strands colored $i$ and $j$ in a given colored tangle, let $\min(\chi):=\min(i,j)$ denote the minimum of the colors involved in the crossing $\chi$.
\end{definition}

With Definition \ref{def-crossingMin} in hand, we have the following handy formula that computes all homological shifts that we will need in this paper.  Note that, as in Remarks \ref{rmk-0edges} and \ref{rmk-0edgesAgain}, this formula still holds when we allow edges labelled by a $0$ to stand in for the lack of an actual edge.

\begin{proposition}\label{prop-isotopyShifts}
Let $T_1$ and $T_2$ be any two colored positive tangles such that we can transform $T_1$ into $T_2$ via a sequence of braid-like Reidemeister moves and fork/ladder twists and slides.  Then $C_N(T_1)\simeq t^\alpha q^\beta C_N(T_2)$ where
\begin{equation}\label{eq-isotopyShifts}
\alpha = \sum_{\chi\in T_1} \min(\chi) - \sum_{\delta\in T_2} \min(\delta).
\end{equation}
Here each sum is taken over all crossings in each diagram.
\end{proposition}
\begin{proof}
All of the shifts for Reidemeister moves and fork/ladder twists and slides indicated throughout this section can be seen to follow precisely this format.
\end{proof}

Equation \ref{eq-isotopyShifts} will often be the simplest way to understand what sort of shift has occurred after performing an isotopy of a given tangle.  Often after such an isotopy, all of the crossings in the diagram will remain in the same place, but some will have `changed color', and this will immediately describe the resulting homological shift $\alpha$.  One such example is illustrated below.

\begin{figure}[ht]
\centering
\includegraphics[scale=.35]{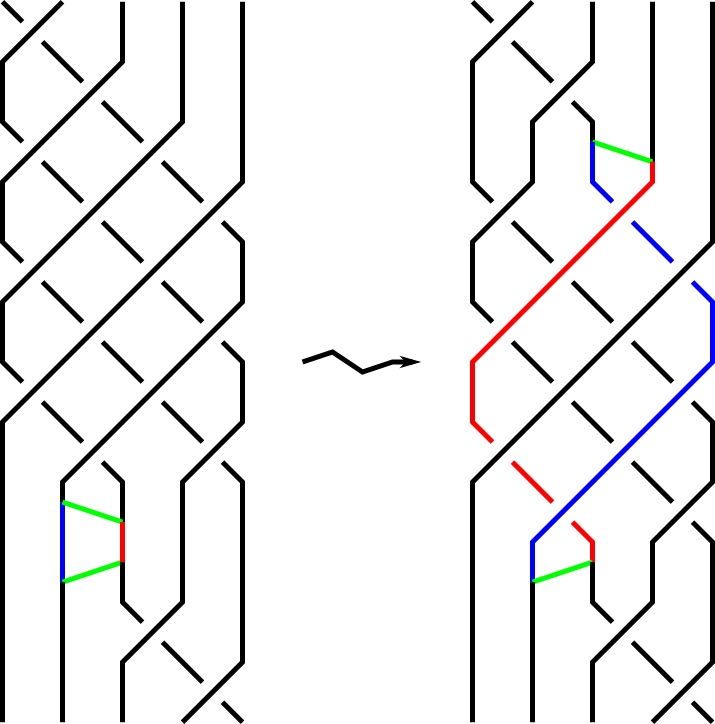}
\caption{The upper rung of the ladder is `pulled through' the twisting via ladder twists and slides, giving a braid-like isotopy between the two colored tangles that leaves all crossings where they are, but changes some colorings. Here the black strands are all the same color, but the blue and green strands are colored strictly less than the black, while the red strand is strictly greater than the black.}
\label{fig-pullRungEx}
\end{figure}

In Figure \ref{fig-pullRungEx}, the black strands are all colored $m$, while the green and blue strands are all colored strictly less than $m$ ($i$ and $j:=m-i$, respectively).  Meanwhile, the red strands are colored $m+i$.  Thus any crossing that maintained its color, and any crossing that became black-red, did not change its min value.  However, crossings that changed to black-blue or red-blue changed by precisely $i$.  Thus in this example with five such crossings, $\alpha=5i$.  We will use similar conventions for the colors of strands (black, red, green and blue) in future diagrams.

\subsection{Homological algebra}
We now recall some key facts and lemmas about the homological algebra of limits and inverse systems from \cite{Roz10a} and \cite{IW16}. For this section, we will assume all chain complexes are over some additive category (say a chain complex of webs in $\foamC$). As with the complexes given above, we will assume our differential increases homological degree by 1.

\begin{definition}\label{def-homOrder}
Let $\AS$ and $\BS$	be chain complexes and suppose $f: \AS\to \BS$ is a chain map. Define the \emph{homological order} of $f$, which we denote by $|f|_h$, to be the maximal degree $d$ for which $[\AS\xrightarrow{f} \BS]$ is chain homotopy equivalent to a complex $\AS'$ that is contractible below homological degree $d$.
\end{definition}

Roughly speaking, $|f|_h$ denotes the maximal homological degree through which $f$ gives a homotopy equivalence on the truncated complexes of $\AS$ and $\BS$ up to degree $|f|_h$.

\begin{definition}\label{def-invSys}
	An \emph{inverse system} of chain complexes is a sequence of chain complexes equipped with chain maps $$\{\AS_k,f_k\}=\AS_0 \xleftarrow{f_0} \AS_1 \xleftarrow{f_1} \AS_2 \xleftarrow{f_2} \cdots.$$ An inverse system is called \emph{Cauchy} if $|f_k|_h \to \infty$ as $k \to \infty$.
\end{definition}

\begin{definition} \label{def-invLim}
An inverse system $\{\AS_k, f_k\}$ has a \emph{limit} (or \emph{inverse limit}), which we denote by $\AS_\infty$ or $\lim \AS_k$, if there exist maps $\tilde{f}_k: \AS_\infty \to \AS_k$ such that
\begin{enumerate}
 	\item The $f_k\circ \tilde{f}_{k+1}=\tilde{f}_k$ for all $k\geq 0$.
 	\item $|\tilde{f}_k|_h \to \infty$ as $k \to \infty$. 
\end{enumerate}
\end{definition}

\begin{theorem}[Rozansky \cite{Roz10a}]\label{thm-Cauchy}
An inverse system $\{\AS_k,f_k\}$ of chain complexes has a limit $\AS_\infty$ if and only if it is Cauchy.
\end{theorem}

We conclude this section with a result from \cite{IW16} which allows us to prove two inverse systems have equivalent limits. 

\begin{proposition}[\cite{IW16}] \label{prop-compSys}
	Suppose $\{\AS_k,f_k\}$ and $\{\BS_{\ell},g_{\ell}\}$ are two Cauchy inverse systems with limits $\AS_\infty$ and $\BS_\infty$ respectively. Let $z(\ell)$ be a nondecreasing function on $\NN_0$ such that $\lim_{\ell\rightarrow\infty}z(\ell)=\infty$. Suppose there are maps $F_\ell : \BS_\ell \to \AS_{z(\ell)}$ forming a commuting diagram with the system maps $f_k$ and $g_{\ell}$ for appropriate $k$ and $\ell$. If $|F_\ell|_h \to \infty$ as $\ell \to \infty$ then $\AS_\infty \simeq \BS_\infty$.
\end{proposition}

\begin{figure}[ht]
\centering
\includegraphics[scale=.5]{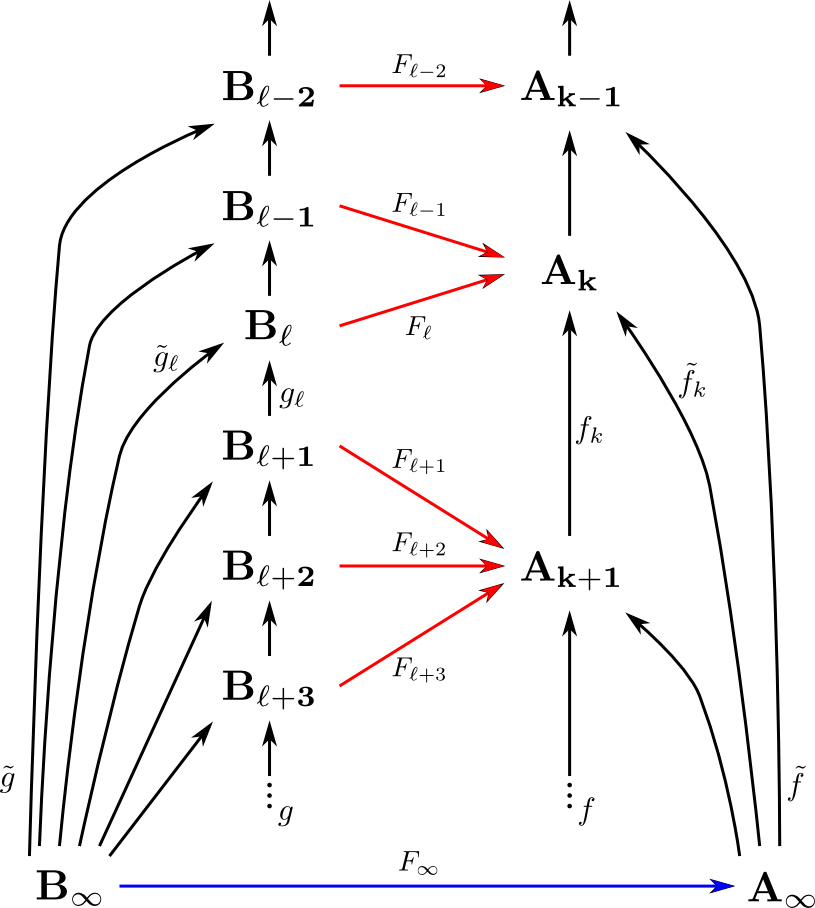}
\caption{The diagram for Proposition \ref{prop-compSys}. In this diagram both $\{\mathbf{A_k},f_k\}$ and $\{\mathbf{B_{\ell}},g_{\ell}\}$ are Cauchy inverse systems. The limits $\mathbf{A_\infty}$ and $\mathbf{B_\infty}$ and the maps $\tilde{f},\tilde{f}_k,\tilde{g},\tilde{g_{\ell}}$ all exist from Theorem \ref{thm-Cauchy}. If we can find the maps $F_\ell$, then Theorem \ref{thm-Cauchy} also provides the map $F_\infty$. We simply need to show that $|F_\ell|_h \to \infty$ as $\ell\to\infty$ to prove that $F_\infty$ is a homotopy equivalence.  \label{fig-invSys}}	
\end{figure}

In Figure \ref{fig-invSys} we include a diagram to better explain the situation in Proposition \ref{prop-compSys}.  In principal we could allow cases where $z(\ell+1)>z(\ell)+1$, necessitating the use of multiple system maps $f_k$ in the commutation required; in the cases of interest in this paper, however, this will never be the case and $z$ will only ever increase by one or not at all.

\subsection{Infinite full twists} \label{subsec-infTwists}

Ultimately, we will use Proposition \ref{prop-compSys} to compare our Cauchy inverse systems to a well understood Cauchy inverse system coming from studying the infinite full twist which we now discuss. 

\begin{figure}[ht]
\centering
\includegraphics[scale=.28]{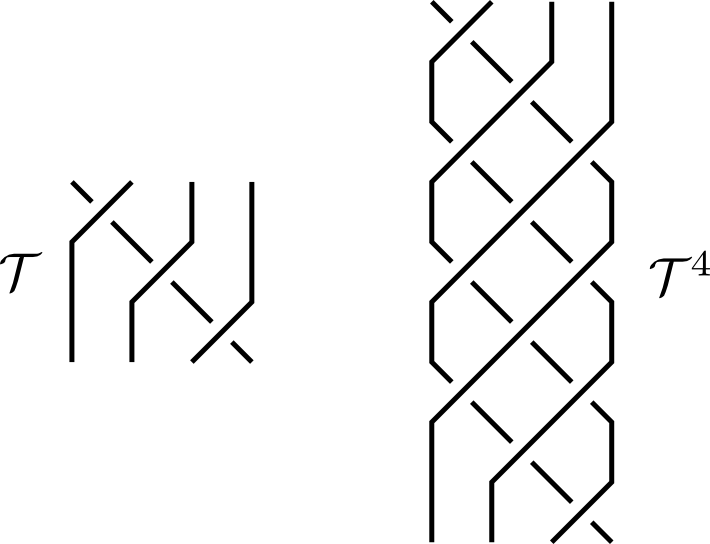}
\caption{Examples of $\T$ and $\FT=\T^4$ for $n=4$ \label{fig-FT}}
\end{figure}

Let $n$ be a fixed positive integer. Let $\Bn$ denote the braid group on $n$ strands and $\sigma_1,...,\sigma_{n-1}$ denote its standard (positive) generators. We let $\T$ denote the fractional twist braid $\sigma_{1}\cdots\sigma_{n-1}$, so that $\FT:=\T^n$ is the full-twist braid on $n$ strands (see Figure \ref{fig-FT} for a specific example). 

For the reminder of this section, we will consider braids $\beta \in \Bn$ which are colored by a single label $m$. That is, all strands of $\beta$ are labeled by the same integer $m \in \{1,...,N-1\}$. We will use the notation $\beta_{(m)}$ to denote such a braid. 

In \cite{Roz10a}, Rozansky proved that the Khovanov complex of the ``infinite full twist'' was a well-defined complex and that it categorified the Jones-Wenzl projector.  We now state the analogous result in colored Khovanov-Rozansky homology due to Cautis.

\begin{theorem}[Cautis \cite{Cau12} ]\label{thm-infFT}   
	Let $\FT_{(m)} \in \Bn$ denote the positive full twist braid on $n$ strands with all strands labeled by $m$, then $$\C{\FT_{(m)}^\infty} := \lim_{k \to \infty} \C{\FT^k_{(m)}}$$ is a well-defined chain complex. Furthermore, $\C{\FT_{(m)}^\infty}$ is an idempotent chain complex categorifying a highest weight projector $p_{n,(m)}: (\bigwedge^m\CC_q^N)^{\otimes n}\to(\bigwedge^m\CC_q^N)^{\otimes n} $ factoring through the unique highest weight subrepresentation of $(\bigwedge^m\CC_q^N)^{\otimes n}$. 
\end{theorem}

We will see in the sequel that we may relax the requirement that the constituent complexes of the inverse system above are full twists. In particular, 
$$\C{\T_{(m)}^\infty} := \lim_{k \to \infty} \C{\T^k_{(m)}} \cong \C{\FT_{(m)}^\infty}.$$ In this text, we will write $P_{n,(m)}:= \C{\FT_{(m)}^\infty}$ for the sake of brevity.

\begin{remark} \label{rmk-Rickard}
Strictly speaking, Theorem \ref{thm-infFT} is proven for Rickard complexes in the categorification of $\qsln$. However this category is equivalent to $\foamC$ and thus can be viewed as a result for complexes associated to full twist braids in $\foamC$ (See \cite{Cau12} and \cite{QR16} for more details).
\end{remark} 

Finally, we may also color the strands of $\FT$ with any set of labels, and the above theorem holds (as proven in \cite{Cau12}). However, we focus on the case of ``unicolored'' braids in this text.

The action of tensoring the complex for a braid generator with the categorified highest weight projector $P_{n,(m)}$ has a very simple description. We will use the following lemma (Proposition 5.8 from \cite{Cau12}) to prove certain corollaries in \S \ref{subsec-mainResult}. 

\begin{lemma} \label{lem-ladKilling}
	$\sladderC \otimes P_{n,(m)} \simeq  P_{n,(m)} \otimes \sladderC \simeq 0.$ In particular, $$C_N(\sigma_{i,(m)}) \otimes P_{n,(m)} \simeq P_{n,(m)}\text{  and  } C_N(\sigma^{-1}_{i,(m)}) \otimes C_N(\FT^\infty_{(m)}) \simeq t^m q^m P_{n,(m)}.$$
\end{lemma}

\section{Main results} In this section, we restate one of our main results and also give some corollaries which allow our main theorem to be applied to other similar cases. After stating these results, we outline the proof of the main theorem which will hinge upon Proposition \ref{prop-compSys} and Figure \ref{fig-invSys}.  The technical brunt of the proof will be handled via several lemmas estimating the homological order of maps between the inverse systems of Proposition \ref{prop-compSys} (the red maps $F_\ell$ of Figure \ref{fig-invSys}).

\subsection{Main theorem and corollaries} \label{subsec-mainResult} We first recall the definition of our main class of braids. 

\begin{definition} \label{def-semiInfBraid}
	A \emph{semi-infinite positive braid word} $\B$ on $n$ strands is a semi-infinite word on the generators $\sigma_i$ of the braid group $\Bn$. 
	$$\B = \sigma_{i_1}\sigma_{i_2}\sigma_{i_3}\cdots$$
	Such a word $\B$ represents a \emph{semi-infinite positive braid} $\BB$ up to allowable braid moves (ie, up to isotopies and Reidemeister moves that do not introduce negative crossings).  We call a semi-infinite positive braid word $\B$ \emph{complete} if every generator $\sigma_i$ of $\Bn$,  for $i=1,\dots,n-1$, appears an infinite number of times in $\B$.
\end{definition}

Note that the property of completeness is preserved by braid relations for positive braids, and so it makes sense to refer to a complete semi-infinite positive braid.

\begin{definition} \label{def-partBraid}
Let $\B = \sigma_{i_1}\sigma_{i_2}\sigma_{i_3}\cdots$ be a semi-infinite positive braid word and $\ell \in \NN$. We define the \emph{$\ell$th partial braid} of $\B$, denoted by $\B_\ell$, as $$\B_\ell = \sigma_{i_1}\sigma_{i_2}\cdots\sigma_{i_\ell}.$$
\end{definition}

In this text we are also concerned with colored braids. Following the notation in \S \ref{subsec-infTwists}, we use the notation $\BB_{(m)}$ to denote the semi-infinite positive braid $\BB$ with every strand labeled by $m$.  Representing $\BB$ by the word $\B$, we use the notation $\B_{\ell, (m)}$ to denote the $\ell$th partial braid with every strand colored by $m$.  We will write $\C{\B_{(m)}}:= \lim_{\ell \to \infty} \C{\B_{\ell,(m)}}$ if the limit exists.  Finally, if the limit $\C{\B_{(m)}}$ is independent of the choice of the representative word $\B$, we can define the complex $\C{\BB_{(m)}}:=\C{\B_{(m)}}$.  We now state the main result of our paper.

\begin{theorem}\label{thm-infBraid}
	Let $\B$ be a complete semi-infinite positive braid word on $n$ strands.  Then $\C{\B_{(m)}}$ exists for all $m \in \{1,...,N\}$. Furthermore, $$\C{\B_{(m)}} \simeq P_{n,(m)}.$$
	Thus, for any complete semi-infinite positive braid $\BB$, we have a well-defined
	$$\C{\BB_{(m)}} \simeq P_{n,(m)}.$$
\end{theorem}

Before proving Theorem \ref{thm-infBraid}, we also record a few corollaries.
  
\begin{corollary}\label{cor-fewNegXs}
Let $\B$ be a complete semi-infinite braid word containing a finite number $a$ of negative crossings. Then
$$\C{{\B}_{(m)}} \simeq t^{ma}q^{ma}P_{n,(m)}$$
and thus, for any complete semi-infinite braid $\BB$ that can be represented by a word having finitely many negative crossings,
 $$\C{{\B}_{(m)}} \simeq P_{n,(m)}$$
up to an overall degree shift.
\end{corollary}
Note that, in the presence of negative crossings, our complexes for finite braids are well-defined only up to these same degree shifts just as well (see Equation \ref{eq-shiftRMoves}).

\begin{proof}
Suppose that the $r$th crossing of $\B$ is the $a$th (and final) negative crossing of $\B$. Then we can write $\B$ as a product $\B_r \B'$ where $\B' = \sigma_{i_{r+1}}\sigma_{i_{r+2}}\sigma_{i_{r+3}}\cdots$ is a positive semi-infinite complete braid. Therefore $\C{\B'_{(m)}} \simeq P_{n,(m)}$ by Theorem \ref{thm-infBraid}. Since $\B_r$ contains exactly $a$ negative crossings and some number of positive crossings, Lemma \ref{lem-ladKilling} implies $$\C{\B_{(m)}}\simeq\C{\B_{r,(m)}}\otimes P_{n,(m)} \simeq (tq)^{am}P_{n,(m)}$$ as desired.
\end{proof}

We can also consider bi-infinite positive braid words of the form
\begin{equation}\label{eq-biinf}
 	\B = \cdots \sigma_{j_{-2}}\sigma_{j_{-1}}\sigma_{j_{1}}\sigma_{j_2}\cdots
\end{equation}
by viewing $\B$ as a product of semi-infinite braid words, say $\B^-\B^+$, where $\B^+ = \sigma_{j_{1}}\sigma_{j_2}\cdots$ and $\B^- = \cdots \sigma_{j_{-2}}\sigma_{j_{-1}}$. In this case we define $$\C{\B_{(m)}} :=\C{\B^-_{(m)}}\otimes\C{\B^+_{(m)}},$$ in the case that both $\C{\B^-_{(m)}}$ and $\C{\B^+_{(m)}}$ exist.

\begin{corollary}\label{cor-biinf}
Let $\BB$ be a complete positive bi-infinite braid represented by a word $\B$ of the form (\ref{eq-biinf}). Then we have a well-defined $\C{\BB_{(m)}} \simeq P_{n,(m)}$ independently of the choice of $\B$.
\end{corollary}

\begin{proof}
Theorem \ref{thm-infBraid} proves that, for any choice of braid word $\B^-\B^+$ representing $\BB$, $\C{\B^+_{(m)}} \simeq P_{n,(m)}$. Likewise, via a similar argument, $\C{\B^-_{(m)}} \simeq P_{n,(m)}$. The result follows from the idempotency of $P_{n,(m)}$.
\end{proof}

Finally we can also study infinite braids $\BB$ which are not complete. We do this by studying complete subbraids of $\BB$ on a subset of the strands. To do so we need to first introduce some notation. 

Let $\B$ be a semi-infinite positive braid word. If $\B$ is not complete, then let 
$$\gamma(\B)=\{i|\sigma_i\text{ appears infinitely often in } \B\} \subset \{ 1,...,n-1\}$$
 and $\gamma^c(\B) = \{0,1,...,n-1,n\}\backslash \gamma(\B)$. 

\begin{definition}\label{def-sqcup}
Let $\BB^1$ be a (possibly bi-infinite or semi-infinite) positive braid on $n_1$ strands and $\BB^2$ be a (possibly bi-infinite or semi-infinite) positive braid on $n_2$ strands. The braid $\BB^1 \sqcup \BB^2$ on $n_1 + n_2$ strands is the braid defined by taking $\BB^1$ on the first $n_1$ strands and $\BB^2$ on the last $n_2$ strands.
\end{definition}

Recall the monoidal product $\sqcup$ defined on $\webC$ and inherited by $\foamC$. It can be directly verified that $\C{\BB^1 \sqcup \BB^2} \simeq \C{\BB^1} \sqcup \C{\BB^2}$, and similarly for colored braids.

\begin{corollary}\label{cor-nonComplete}
	Let $\BB$ be a semi-infinite positive braid represented by $\B$, and suppose $\gamma^c(\B) = \{i_0 =0,i_1,...,i_{p-1}, i_p = n\}$. Let $\B_r$ be the partial braid of $\B$ containing all copies of the generators $\sigma_{i_1},...,\sigma_{i_p}$, so that none of these generators appear after the $r$th crossing of $\B$. Then $\C{\B_{(m)}} \simeq (P_{n_1,(m)} \sqcup \cdots \sqcup P_{n_p,(m)}) \otimes \C{\B_{r,(m)}}$ giving a well-defined
	$$\C{\BB_{(m)}} \simeq (P_{n_1,(m)} \sqcup \cdots \sqcup P_{n_p,(m)}) \otimes \C{\B_{r,(m)}}.$$
	Here $n_j = i_j -i_{j-1}$ for $j = 1,...,p$ and any terms with $n_j = 0$ are omitted.
\end{corollary}

\begin{figure}[ht]
\includegraphics[scale=.3]{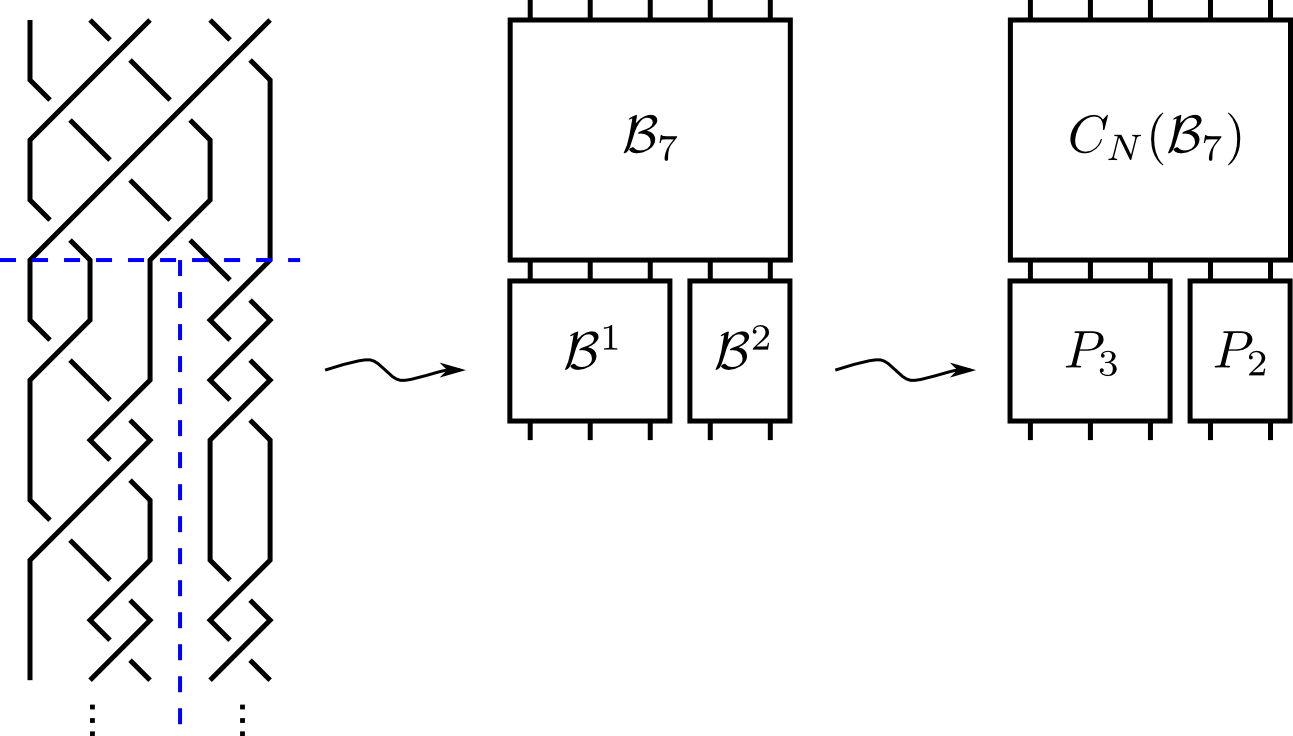}
\caption{An example of a non-complete semi-infinite braid $\BB$ is pictured on the left.  The dashed lines indicate how we consider partitioning this braid into the partial braid $\B_7$ and the two complete semi-infinite braid words $\B^1$ and $\B^2$, leading to projectors $P_3$ and $P_2$ respectively.  In the notation of Corollary \ref{cor-nonComplete}, $\gamma^c(\B)=\{0,3,5\}$ and $r=7$ indicates the final appearance of the braid group generator $\sigma_3$ after indexing as in the partitions.  Colors are omitted.}
\label{fig-nonComplete}
\end{figure}

The proof of this corollary follows from realizing $\B = \B'\B_r$ where $\B'$ can be written as $\B^1 \sqcup \cdots \sqcup \B^p$ for some semi-infinite complete positive braid words and $\B_r$.  All of this is clearly upheld by braid relations, and thus we have well-defined complexes.  See Figure \ref{fig-nonComplete} for a simple example.

\subsection{Proof of the main theorem}\label{subsec-proofMain}

To prove Theorem \ref{thm-infBraid}, we will be following an outline similar to that of the proof of the main result in \cite{IW16}. More precisely, we will use a combination of Proposition \ref{prop-isotopyShifts} and Proposition \ref{prop-compSys} with the following correspondences to Figure \ref{fig-invSys} in mind.

\begin{enumerate}
\item The complexes $\C{\FT_{(m)}^{k}}$ will play the role of the $\mathbf{A_k}$.
\item Therefore by Theorem \ref{thm-infFT}, $\mathbf{A_\infty} \simeq P_{n,(m)}$.
\item Given a complete semi-infinite positive braid word $\B$, the complex $\C{\B_{\ell,(m)}}$ plays the role of the $\mathbf{B_\ell}$.
\item The maps $g_\ell$ will be precisely a tensor product of identity maps with the canonical projection $\pi$ of Equation \ref{eq-posxProj} (see Lemma \ref{lem-posxCone}) corresponding to the $\ell$th (and final) crossing in $\B_\ell$.  Thus $g_\ell$ corresponds to resolving the $\ell$th crossing in $\B_\ell$, bringing us to $\B_{\ell-1}$.
\item The maps $F_\ell$ will be constructed in a similar manner by resolving certain carefully chosen crossings in $\B_\ell$.
\item The function $k = z(\ell)$ will be based upon how many full twists $\FT^k$ can be found as a subbraid of $\B_\ell$. The condition that $\B$ is complete guarantees that $k\rightarrow\infty$ as $\ell\rightarrow\infty$.
\item The estimates on $|F_\ell|_h$ will be based on Proposition \ref{prop-isotopyShifts} (and specifically Equation \ref{eq-isotopyShifts}) and Lemma \ref{lem-posxCone}. We will show the system $\{\mathbf{B_\ell},g_\ell\}$ is Cauchy in a similar manner.
\end{enumerate}

We now proceed to make these points more precise.  The first four points need no further elaboration.  We begin by handling the fifth and sixth points above, that is, the construction of the map $F_\ell$ and the corresponding function $k(\ell)$ determining the target of $F_\ell$.  This portion of the proof is nearly identical to the construction in \cite{IW16}.

\begin{figure}
\begin{center}
\includegraphics[scale=.42]{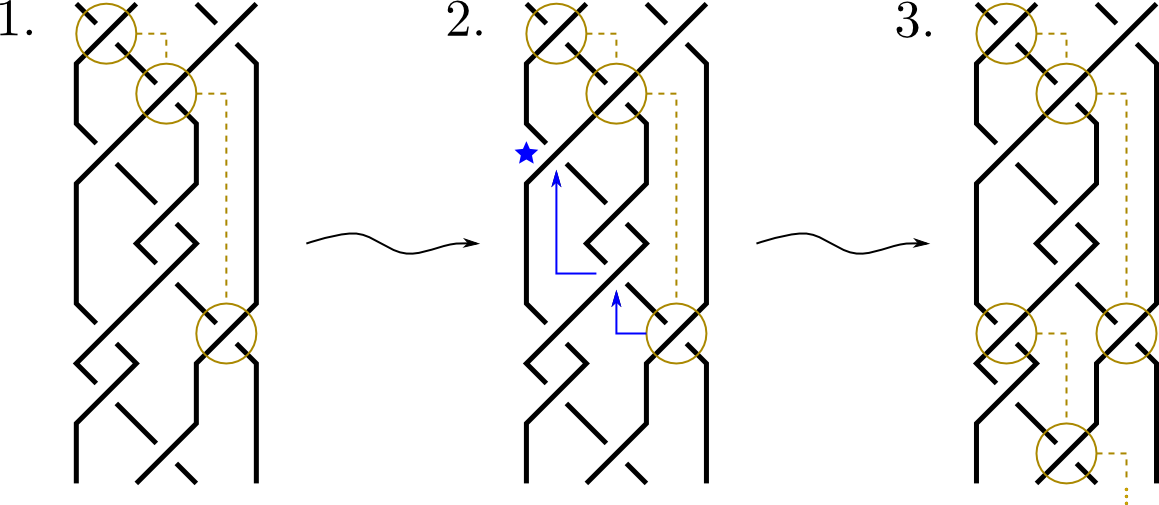}
\end{center}
\caption{An illustration of finding diagonals within some $\B_\ell$.  In step 1, we find the first diagonal illustrated in gold.  In step 2, we work our way back up from the diagonal $\sigma_{n-1}$ as in the blue arrows until we arrive at the $\sigma_1$ marked by a blue star.  In step 3, we begin forming the second diagonal starting from the first $\sigma_1$ below the starred crossing from step 2. In this way all of the gold dashed lines are kept disjoint from one another.}
\label{creating diagonals}
\end{figure}

Given the partial braid $\B_\ell$, we start at the top of the braid (beginning of the braid word) and seek the first occurrence of the generator $\sigma_1$.  From that point we go downward and find the first occurrence of $\sigma_2$, and so forth until we reach $\sigma_{n-1}$.  In this way we have found crossings within $\B_\ell$ that would, in the absence of the crossings we `skipped', give a single copy of $\T^1$.  We connect these crossings with a dashed line going rightward then downward as in Figure \ref{creating diagonals}, and we call the crossings involved \emph{diagonal crossings}.  Having found such a diagonal within $\B_\ell$, we work our way back \emph{up} the braid $\B_\ell$ in the same way going from the diagonal $\sigma_{n-1}$ to the previous (not necessarily diagonal) $\sigma_{n-2}$ and so forth until we reach another $\sigma_1$ (if there were no skipped crossings, we are now back at the $\sigma_1$ we started with).  We begin the second diagonal from the first $\sigma_1$ that is \emph{below} this $\sigma_1$ we found at the end of our upward journey.  In this way we find disjoint diagonals with as few `skipped' crossings between them as possible.  See Figure \ref{creating diagonals} for clarification.

Let $y(\ell)$ denote the number of diagonals that can be completed within $\B_\ell$ in this way.  The function $z(\ell)$ determining the destination of the map $F_\ell$ is
\begin{equation}
\label{z(ell) eqn}
z(\ell):=\left\lfloor\frac{y(\ell)}{n}\right\rfloor
\end{equation}
where $\lfloor\cdot\rfloor$ denotes the integer floor function.  Thus $z(\ell)$ gives the number of \emph{full} twists that can be seen within $\B_\ell$.

At this point, we `erase' the diagonals below the $nz(\ell)^{\text{th}}$ diagonal - that is, all crossings below the last diagonal of the last full twist within $\B_\ell$ are regarded as non-diagonal.  The map $F_\ell$ is then the composition of mapping cone projections $\pi$ coming from Lemma \ref{lem-posxCone} where we are resolving all non-diagonal crossings in $\B_\ell$.  Note that the order in which we resolve the crossings is irrelevant, however for the sake of consistency here we will assume that we resolve crossings starting from the bottom upwards.  From this consideration it should be clear that the maps $F_\ell$ commute with the maps $f_k$ and $g_\ell$ of the two systems $\C{\FT_{(m)}^k}$ and $\C{\B_{\ell,(m)}}$, which are also just maps based on resolving bottom-most (final) crossings.

\begin{figure}
\begin{center}
\includegraphics[scale=.42]{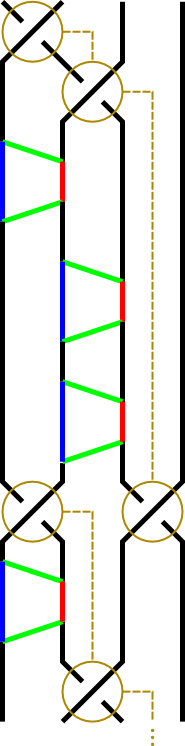}
\end{center}
\caption{An example of a possible resolution of the diagram from Figure \ref{creating diagonals} involving various ladders between the diagonals.  In this diagram, the black strands are all colored $m$.  Similar to the coloring in Figure \ref{fig-pullRungEx}, the blue and green strands are all colored strictly less than $m$, while the red strands are colored strictly greater than $m$.  Although any two blue strands may actually have different labels depending on the resolution, they both must be less than $m$ (and similarly for the greens and reds).}
\label{fig-exampleResolution}
\end{figure}

For the seventh and final point listed above, we seek to estimate $|F_\ell|_h$ for the $F_\ell$ constructed above.  To do this, we consider the mapping cone of one projection $\pi$ (resolving one crossing) at a time.  As discussed below Lemma \ref{lem-posxCone}, this mapping cone consists of a complex involving ladders in place of the resolved crossing.  We then expand all of the terms in this complex by resolving the other non-diagonal crossings.  In this way, we exhibit $\cone(\pi)$ as the total complex of a large number of chain maps where each term consists of the chain complex of a diagram involving the diagonal crossings together with \emph{at least one} ladder present in place of one of the non-diagonal crossings (every one of these terms but for the `all-straight' resolution of the remaining crossings will contain several ladders). Figure \ref{fig-exampleResolution} provides an example of such a diagram that could result from resolving the non-diagonal crossings in Figure \ref{creating diagonals}.  Estimating $|\pi|_h$ (and thus $|F_\ell|_h$ which comes from a composition of projections $\pi$) is the same thing as estimating the minimum homological degree of any such complex.

In order to find a lower bound on this minimum homological degree, we first introduce some further terminology matching the terminology used in \cite{IW16}.

\begin{definition}\label{def-diagsAndZones}
Given a partial braid $\B_{\ell,(m)}$ coming from some uni-colored complete semi-infinite positive braid word $\B_{(m)}$, let $D$ be any colored diagram arrived at from $\B_{\ell,(m)}$ by resolving the non-diagonal crossings.  The dashed gold lines between diagonal crossings (as in Figure \ref{creating diagonals}) are called \emph{diagonals}.  The space between any two such diagonals will be called a \emph{zone}.  A zone is called \emph{empty} if there are no ladders contained within it; otherwise, it is called \emph{non-empty}.  In any non-empty zone, the \emph{top-most} (respectively \emph{bottom-most}) ladder indicates the ladder replacing the non-diagonal crossing $\sigma_{i_j}$ of minimal (respectively maximal) $j$ resolved within the given zone.
\end{definition}
The diagram with all empty zones is precisely the diagram $\FT^{k(\ell)}$ that is used to determine the target of $F_\ell$ (recall that we resolve all crossings below the final diagonal of the final full twist within $\B_\ell$).  As discussed above, estimating $|F_\ell|_h$ is the same as bounding the minimum homological degree of any complex involving only diagrams with at least one non-empty zone.

At this point, our proof diverges slightly from that of \cite{IW16} and we need a more careful, more conservative estimate on these minimum homological degrees.  To begin with, we note that any ladder in any of the diagrams $D$ involves one vertical edge with label greater than $m$, and one vertical edge with label less than $m$ (for simplicity, we draw the vertical edge labelled zero if necessary).  Thus when we perform braid-like isotopies as in Proposition \ref{prop-isotopyShifts}, we will always be in situations similar to that of Figure \ref{fig-pullRungEx}.

\begin{figure}[ht]
\begin{center}
\includegraphics[scale=.42]{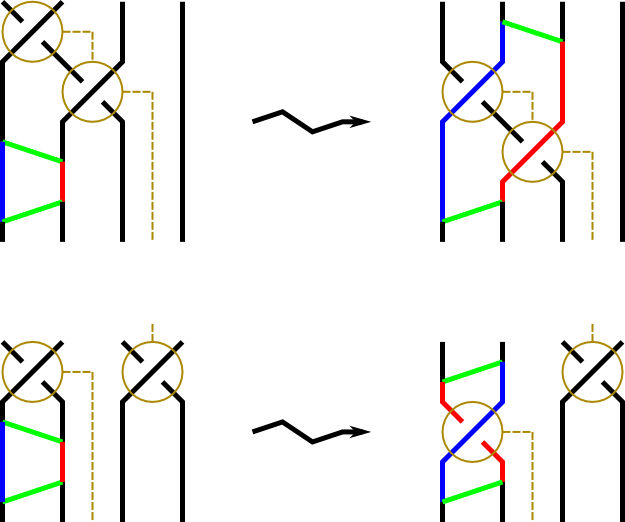}
\end{center}
\caption{The two cases of Lemma \ref{lem-nonemptyZonesBound}.  The first case is illustrated via the first non-empty zone from the diagram of Figure \ref{fig-exampleResolution}, where a single ladder slide is available.  The second case is illustrated via the second zone from Figure \ref{fig-exampleResolution} involving a single ladder twist.}
\label{fig-UpperLadderMove}
\end{figure}

\begin{lemma}\label{lem-nonemptyZonesBound}
For any colored diagram $D$ arrived at from $\B_{\ell,(m)}$ as in Definition \ref{def-diagsAndZones}, let $b_1(D)$ indicate the number of non-empty zones in $D$, omitting the zone above the first diagonal.  Then the minimum homological degree $\min_h(\C{D})$ has lower bound
\begin{equation}\label{eq-nonemptyZonesBound}
\min_h(\C{D})\geq b_1(D).
\end{equation}
\end{lemma}
\begin{proof}
The upper rung of the top-most ladder in every non-empty zone (except for the top one) can be pulled through the diagonal above it, either with a ladder slide or with a ladder twist.  Both cases are illustrated in Figure \ref{fig-UpperLadderMove} which uses the two zones present in the diagram of Figure \ref{fig-exampleResolution}.  In either case we have one (diagonal) crossing with a new, lower minimum color so that Proposition \ref{prop-isotopyShifts} guarantees a shift of at least 1 in homological degree.
\end{proof}

\begin{figure}[ht]
\begin{center}
\includegraphics[scale=.42]{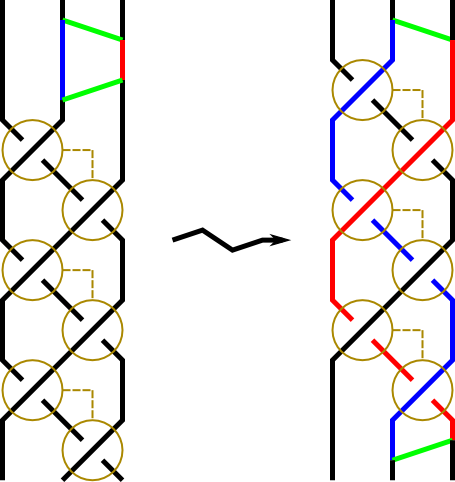}
\end{center}
\caption{Pulling a ladder rung downward through a full twist, as in Lemma \ref{lem-downwardFTBound}.  The case $n=3$ is illustrated here, with the ladder starting above the second and third strands.  All such cases start with a ladder slide similarly, except that if the ladder starts above the first and second strands, then the pull begins with a ladder twist.}
\label{fig-LowerLadderFT}
\end{figure}

\begin{lemma}\label{lem-downwardFTBound}
For any colored diagram $D$ arrived at from $\B_{\ell,(m)}$ as in Definition \ref{def-diagsAndZones}, let $b_2(D)$ indicate the sum of the number of `empty-zoned full twists' below each non-empty zone.  That is, starting from each non-empty zone, we count the number of diagonals $d$ below it that have only empty zones between them, and take $\left\lfloor\frac{d}{n}\right\rfloor$.  Then $b_2(D)$ is the sum of all of these counts for each non-empty zone.  With this, the minimum homological degree $\min_h(\C{D})$ has lower bound
\begin{equation}\label{eq-downwardFTBound}
\min_h(\C{D})\geq 2(n-1)b_2(D).
\end{equation}
\end{lemma}
\begin{proof}
The bottom rung of the bottom-most ladder in a non-empty zone can be pulled through a full twist below it as illustrated in Figure \ref{fig-LowerLadderFT} for $n=3$ strands.  When this is done, the blue edge with the lower color crosses each other strand twice during a full twist, and there are $(n-1)$ other strands (including the red strand with higher color), so that Proposition \ref{prop-isotopyShifts} guarantees a shift of at least 1 per crossing involving the lower color edge.
\end{proof}

\begin{figure}[ht]
\begin{center}
\includegraphics[scale=.42]{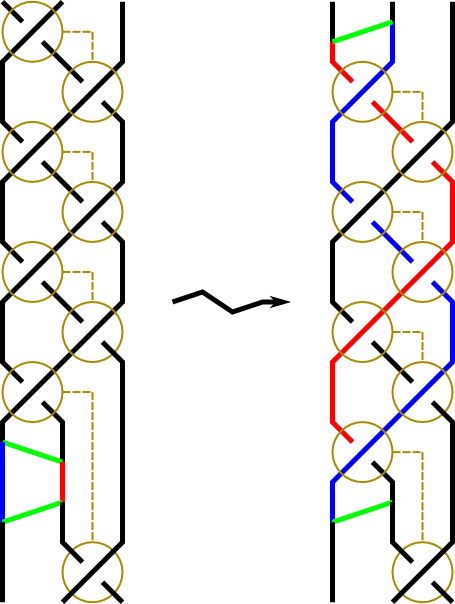}
\end{center}
\caption{Pulling a ladder rung upward through a full twist, as in Lemma \ref{lem-upwardFTBound}.  The diagram illustrates how, due to the placement of the ladder, a ladder twist is needed to pass through the first diagonal before one can count diagonals for full twists.  Note that if the top diagonal were missing from this diagram, the rung would not end up in the same position and the count of crossing color changes would no longer be correct.}
\label{fig-UpperLadderFT}
\end{figure}

\begin{lemma}\label{lem-upwardFTBound}
For any colored diagram $D$ arrived at from $\B_{\ell,(m)}$ as in Definition \ref{def-diagsAndZones}, we start from each non-empty zone and count the number of diagonals $d'$ above it that have only empty zones between them, and take $\left\lfloor\frac{d'-1}{n}\right\rfloor$. Define $b_3(D)$ to be the sum of all of these counts for each non-empty zone.  With this, the minimum homological degree $\min_h(\C{D})$ has lower bound
\begin{equation}\label{eq-upwardFTBound}
\min_h(\C{D})\geq 2(n-1)b_3(D).
\end{equation}
\end{lemma}
\begin{proof}
The proof here is the same using the top rung of the top-most ladder.  However, when pulling rungs upward, if there is a fork twist that must occur before pulling the rung through, we lose one diagonal before beginning our full twists, and hence we count based on $d'-1$ rather than $d'$.  See Figure \ref{fig-UpperLadderFT}.  Once again, Proposition \ref{prop-isotopyShifts} provides the estimate based on this count.
\end{proof}

\begin{figure}[ht]
\begin{center}
\includegraphics[scale=.42]{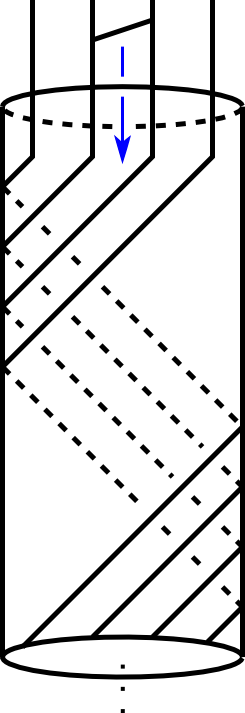}
\end{center}
\caption{An illustration of Remark \ref{rmk-pullingThroughTorus}, where we view the rung of a ladder as being pulled through the interior of the torus defining the full twist as a torus braid.}
\label{fig-pullingThroughTorus}
\end{figure}

\begin{remark}\label{rmk-pullingThroughTorus}
The existence of the braid-like isotopy that pulls a rung through a full twist is perhaps most clearly seen by viewing the rung as passing through the interior of the torus defining a torus braid.  See Figure \ref{fig-pullingThroughTorus}.  In the awkward case of Figure \ref{fig-UpperLadderFT}, the ladder twist needs to be handled first before the rung is to be considered as entering the torus braid.
\end{remark}

\begin{figure}[ht]
\begin{center}
\includegraphics[scale=.42]{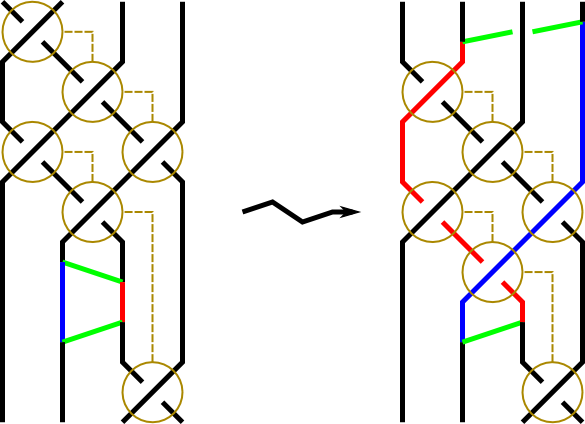}
\end{center}
\caption{An example illustrating the need to count full twists instead of simply counting the number of diagonals as in \cite{IW16}.  Here, if we try to pull the upper-most rung through the two available diagonals, we create a new crossing that did not exist before because of the rung. This new crossing potentially counteracts the shifts coming from the color changes of the diagonal crossings.  It is not difficult to find similar issues when attempting to pull rungs downwards through less than full twists.}
\label{fig-nonFTissue}
\end{figure}

\begin{remark}\label{rmk-newVsOld}
The reader who is familiar with \cite{IW16} may wonder why we have introduced the extra complication of counting our isotopies via full twists of diagonals, rather than simply counting pulls through diagonals which were sufficient in that earlier paper.  The reason for this is that pulling and twisting rungs through anything less than a full twist has the potential to introduce new crossings that were not present before the pull, because of the presence of the rung.  See Figure \ref{fig-nonFTissue}.
\end{remark}

\begin{proof}[Proof of Theorem \ref{thm-infBraid}]
We build a commuting diagram as in Figure \ref{fig-invSys} using the listed points at the beginning of this section.  The construction of $F_\ell$ as a composition of mapping cone projections $\pi$ ensures $|F_\ell|_h$ is at least as large as the minimum $|\pi|_h$ amongst all such $\pi$.  Let $\mathcal{D}_\ell$ denote the set of all diagrams $D$ coming from  $\B_{\ell,(m)}$ as in Definition \ref{def-diagsAndZones} having at least one non-empty zone.  Then, as described above, $|F_\ell|_h$ is bounded below by the minimum homological degree amongst all the $\C{D}$ for all possible diagrams $D\in\mathcal{D}_\ell$.  Lemmas \ref{lem-nonemptyZonesBound}, \ref{lem-downwardFTBound}, and \ref{lem-upwardFTBound} ensure that these homological degrees produce the lower bound
\begin{equation}\label{eq-bigHomBound}
|F_\ell|_h \geq \min_{D\in\mathcal{D}_\ell} \max(b_1(D),2(n-1)b_2(D),2(n-1)b_3(D))
\end{equation}
The assumption that the semi-infinite braid word $\B$ is complete ensures that the number of diagonals in the partial braid $\B_\ell$ increases towards infinity as $\ell\rightarrow\infty$.  The reader may quickly verify that the right hand side of Equation \ref{eq-bigHomBound} must go to infinity as the number of diagonals goes to infinity (roughly, the more diagonals we have, the larger $b_1$ tends to be unless there are many empty zones, in which case $b_2$ or $b_3$ must be large).  Thus we have that $|F_\ell|_h\rightarrow\infty$ as $\ell\rightarrow\infty$.  A similar (and simpler) argument also ensures $|g_\ell|_h\rightarrow\infty$ as $\ell\rightarrow\infty$, verifying that this system is Cauchy and has a limit.  Thus we may use Proposition \ref{prop-compSys} to conclude the proof.
\end{proof}

\section{Infinite braids in HOMFLY-PT homology}
Now we consider our main theorem in the case of HOMFLY-PT homology and its colored generalizations. Up to some technical details which we will discuss, the same methods will be used to prove similar results in these homology theories. 
\subsection{HOMFLY-PT homology}\label{subsec-HHH} First we briefly recall the definition of HOMFLY-PT homology. For the rest of this section, we fix a positive integer $n$. Let $R = \QQ[x_1,...,x_n]$ be a graded polynomial ring with $\deg(x_i) = -2$. Define $R^i$ to be the subring of $R$ which consists of polynomials symmetric in the variables $x_i$ and $x_{i+1}$. We define the Soergel bimodule 
\begin{equation}\label{eq-SoeBi}
 	B_i := q(R \otimes_{R^i} R).
\end{equation}

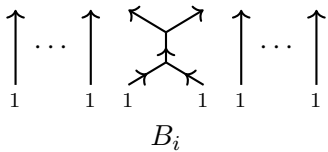
\begin{figure}[ht]
\begin{tikzpicture}
\draw[->,thick] (3,0)--(3,1);
\draw[->-,thick] (1.5,0)--(2,0.3);
\draw[->-,thick] (2.5,0)--(2,0.3);
\draw[->-,thick] (2,.3)--(2,.7);
\draw[->,thick] (2,.7)--(1.5,1);
\draw[->,thick] (2,.7)--(2.5,1);
\draw[->,thick] (1,0)--(1,1);
\draw[->,thick] (4,0)--(4,1);
\draw[->,thick] (0,0)--(0,1);
\node at (0.5,0.5) {$\cdots$};
\node at (3.5,0.5) {$\cdots$};

\node at (2,-0.7) {$B_i$};
\node at (0,-.2) {$_1$};
\node at (1,-.2) {$_1$};
\node at (1.5,-.2) {$_1$};
\node at (2.5,-.2) {$_1$};
\node at (3,-.2) {$_1$};
\node at (4,-.2) {$_1$};
\end{tikzpicture}
\caption{Diagrammatics for $B_i$ in terms of webs. \label{fig-HOMFLYWebs}}
\end{figure}

We can visualize the Soergel bimodules $R$ and $B_i$ in terms of webs. We associate the 1-colored identity web to $R$ and to $B_i$ we associate the web shown in Figure \ref{fig-HOMFLYWebs}. The relations in Figure \ref{fig-webRel} hold the same as in the case of $\sln$-webs (see \cite{Soe92, EKh10} for example). We will use the term ``$\slinf$-webs'' for webs in the case of HOMFLY-PT homology (which we discuss in more detail in \S \ref{subsec-clrHHH}).   

We now define complexes $\F(\b)$ of $R$-$R$-bimodules associated to every braid. These complexes are known as Rouquier complexes. We associate to the positive braid generators $\sigma_i$ the complex
\begin{equation}\label{eq-posxRou}
	\F(\sigma_i) = 0 \xrightarrow{} R \xrightarrow{\chi_i} tq B_i \xrightarrow{} 0, 	
\end{equation}

where $\chi_i(1) = x_i\otimes 1 - 1 \otimes x_{i-1}$. Also we associate to the negative braid generators $\sigma_i^{-1}$ the complex,
\begin{equation}\label{eq-negxRou}
	\F(\sigma_i^{-1}) = 0 \xrightarrow{} t^{-1}q^{-1}B_i \xrightarrow{\chi_o} R \xrightarrow{} 0, 	
\end{equation}

where $\chi_o(f\otimes g) = fg$. We define the complex for a general braid $\B$ via the relation $\F(\b\b') = \F(\b) \otimes_R \F(\b')$. 

Graphically, we denote these complexes in an analogous manner to those in Equation (\ref{eq-posxDef}) in the $1$-colored case.  See Figure \ref{fig-RouCplxs}.

\begin{figure}[ht]
\begin{equation*}
\begin{aligned}
&\posxT = \arcsT \xrightarrow{\chi_i} tq\barbell{}{}{}{}{}\\
&\negxT = (tq)^{-1}\barbell{}{}{}{}{} \xrightarrow{\chi_o} \arcsT\\
\end{aligned}	
\end{equation*}
\caption{Complexes $\F(\sigma_i)$ and $\F(\sigma_i^{-1})$ presented diagrammatically.  \label{fig-RouCplxs}}
\end{figure}

To define a complex for a braid closure we apply the Hochschild homology functor $\HH(\bullet) = H(\bullet \otimes^\text{\textbf{L}}_{R\otimes_\QQ R} R)$ to $\F(\b)$. More precisely, if $\F(\b)$ is the chain complex of bimodules 
$$\cdots \xrightarrow{d} \F_i(\b) \xrightarrow{d} \F_{i+1}(\b) \xrightarrow{d} \cdots,$$
where $\F_i(\b)$ is the component of $\F(\b)$ in homological degree $i$, then $\HH(\F(\b))$ is the complex of vector spaces
$$\cdots \xrightarrow{\HH(d)} \HH(\F_i(\b)) \xrightarrow{\HH(d)} \HH(\F_{i+1}(\b)) \xrightarrow{\HH(d)} \cdots.$$
We then define $\H(\b) := H(\HH(\F(\b)))$. $\H(\b)$ is a triply-graded vector space, with gradings coming from the internal grading on $R$, the homological grading, and a new grading introduced by Hochschild homology which we will denote by $a$ (that is, we will write $a^j C$ to denote that the object $C$ is shifted up in Hochschild degree by $j$).

\begin{theorem}[Khovanov \cite{Kh07}]  \label{thm-HHH}
Let $L$ be a link and let $\b\in \Bn$ be a braid representative for $L$. Then, up to a monomial grading shift depending on the writhe and braid index of $\b$, $\H(\b)$ is a link invariant. Let $\mathcal{P}(\b)$ denote the Poincar\'e series of $\H(\b)$. That is $$\mathcal{P}(\b)(t,q,a) = \sum_{i,j,k\in \ZZ} d(i,j,k)t^iq^ja^k,$$ where $d(i,j,k)$ is the dimension of $\H(\b)$ in grading $(i,j,k)$. Then $\mathcal{P}(\b)(-1,q,a)$ is the HOMFLY-PT polynomial (up to a renormalization).
\end{theorem}

Furthermore, the relations (\ref{eq-forkSlide}) and (\ref{eq-forkTwist}) follow as well in the case that we color all strands by $i = j = k= 1$. More precisely,

\begin{equation}\label{eq-HHHforkSlide}
		\begin{aligned}
			\F(T_1) &\simeq tq\F(T_1')\\
			\F(T_2) &\simeq \F(T_2'),
		\end{aligned}
	\end{equation}
	
where $T_1,T_1',T_2,T_2'$ are the tangles in Figure \ref{fig-forkSlide} in the case that $i=j=k=1$. Also 

\begin{equation}\label{eq-HHHforkTwist}
		\begin{aligned}
			\F(T_3) &\simeq tq^2\F(T_3')\\
			\F(T_4) &\simeq q^{-1}\F(T_4')
		\end{aligned}
\end{equation}

where $T_3,T_3',T_4,T_4'$ are the tangles in Figure \ref{fig-forkTwist} in the case that $i=j=k=1$. Therefore the analogue of Proposition \ref{prop-isotopyShifts} in the $1$-colored case holds for complexes of Soergel bimodules. Likewise, repeating the argument from \S \ref{subsec-proofMain} proves the following theorem.

\begin{theorem}\label{thm-HHHinfBraidConv}
	Let $\FT^k$ denote the $k$th power of the positive full twist braid on $n$ strands. Define $\AS_k = \F(\FT^k)$ and $f_k:\AS_{k+1} \to \AS_k$ to be the canonical projection map (as defined in (\ref{eq-posxCone}) for the $1$-colored case). Then the inverse system $\{\AS_k,f_k\}$ is Cauchy, and thus has a limit $P=\AS_\infty$.
	
	Now let $\BB$ be a positive semi-infinite complete braid on $n$ strands represented by the word $\B$. Define an inverse system by $\BS_\ell = \F(\B_\ell)$ and take $g_\ell$ to be the canonical projection map $\BS_{\ell+1} \to \BS_\ell$. Then $\{\BS_\ell, g_\ell\}$ has a limit $\F(\B)=\BS_\infty$. Furthermore $\BS_\infty \simeq P$.
\end{theorem}

However, the analogue of $P_{n,(1)}$ for HOMFLY-PT homology was defined in a different manner by Hogancamp in \cite{Hog15}. In particular, he used the notion of homotopy limits.\footnote{Hogancamp actually defines his analogue of $P_{n,(1)}$ in terms of homotopy colimits. However the analogous statements hold for homotopy limits.}

\begin{definition}\label{def-holim}
 	Let $\{\AS_k, f_k\}$ be an inverse system of chain complexes.  The \emph{homotopy limit} of $\{\AS_k, f_k\}$ is by definition the mapping cone
 	\[
 	\holim \AS_k := \left[ \prod_{k=0}^\infty \AS_k \xrightarrow{Id - S}\prod_{k=0}^\infty \AS_k \right]
 	\]
 	Where $S:\prod_{k=0}^\infty \AS_k\rightarrow \prod_{k=0}^\infty \AS_k$ is the chain map with components given by the $f_k$.
 \end{definition}
 
\begin{theorem}[Hogancamp \cite{Hog15}]\label{thm-HogProj} Let $\{\AS_k,f_k\}$ be the inverse system defined in Theorem \ref{thm-HHHinfBraidConv}. Then $\tilde{P}=\holim \AS_k$ is an idempotent complex categorifying the one-row Young symmetrizer in the Hecke algebra.
\end{theorem}

$\tilde P$ has been used in recent work of Elias, Hogancamp, and Mellit to compute the HOMFLY-PT homology of torus links explicitly in a combinatorial fashion \cite{EH16, Hog17, Mel17}. Also $\tilde P$ can be used to construct categorifications of one-row colored HOMFLY-PT polynomials \cite{Hog15}. 

With these applications in mind, we wish to understand the relationship between our $P$ and the $\tilde P$ of Hogancamp. In general $\holim\, \AS_k$ and $\lim \AS_k$ are not directly comparable, though in our case this is not an issue. We first recall a definition and then a theorem of Hogancamp \cite{Hog15}.

\begin{definition}\label{def-BSBimod} A \emph{Bott-Samelson bimodule} is a direct sum of bimodules of the form $$M=q^iB_{i_1}\otimes_R \cdots \otimes_R B_{i_\ell}.$$	
\end{definition}

\begin{theorem}\label{thm-HogProjEU}[Hogancamp \cite{Hog15}]
There exists a bounded from below chain complex $\tilde P$ and a map $\varepsilon: \tilde P \to R$ such that
\begin{enumerate}
\item $\tilde P \otimes_R B_i \simeq B_i \otimes_R \tilde P \simeq 0$.
\item $\left[ \tilde P \xrightarrow{\varepsilon} R\right ]$ is a chain complex of bimodules such that $q^iR$ does not appear as a direct summand of any chain bimodule.

\end{enumerate}
	Furthermore, the pair $(\tilde P, \varepsilon)$ is unique up to canonical homotopy equivalence.
\end{theorem}

By construction of our inverse system, we have a map $\varepsilon= \pi: P \to R$ given by the projection from the inverse limit to $\AS_0= \FT^0$. By construction of the complexes $\AS_k = \FT^k$, $R$ can appear only in homological degree $0$ and thus by construction $\cone(\varepsilon)$ satisfies condition (2) of Theorem \ref{thm-HogProj}. To show that $P$ satisfies condition (1) we recall the following lemma from \cite{Hog15} .

\begin{lemma}\label{lem-FT}
Let $C$ be a bounded from below chain complex whose chain bimodules are direct sums of shifts of Bott-Samelson bimodules. If $C$ is supported in homological degrees $\geq \ell$ then $\FT^k \otimes_R C$ is supported in homological degrees $\geq \ell +2k$.
\end{lemma}

Let $\{\CS_k,h_k\}$ be an inverse system of bounded below chain complexes of bimodules. Suppose there exists integers $c_k$ such that $\CS_k$ is supported is homological degrees $\geq c_k$ and $c_k \to \infty$ as $k \to \infty$. Then $|h_k| \to \infty$ as $k \to \infty$ (the system is Cauchy), and $\lim \CS_k$ is contractible. Letting $\CS_k = \FT^k \otimes_R B_i$, $h_k$ be defined as $f_k \otimes 1$, and combining this with Lemma \ref{lem-FT} we see that $P$ satisfies property (1) thus implying the following theorem.

\begin{theorem}\label{thm-equivofP}
$P$ and $\tilde P$ are homotopy equivalent.	
\end{theorem}

Recent work of Mellit \cite{Mel17} builds off of the work of Elias and Hogancamp \cite{EH16, Hog17} to give an explicit computation for the HOMFLY-PT homology of torus knots. By giving sharper bounds on the convergence rate of $|F_\ell|_h$ we present a computation of the HOMFLY-PT homology of a braid positive link up to a certain homological degree. In the following we use $\H(C)$ as shorthand for $H(\HH(C))$, where $C$ is a complex of Soergel bimodules.

Recall a link $L$ is called \emph{braid positive} if it is the closure of a positive braid $\beta$. In this case, $\F(\beta)$ can be viewed as $\F(\B_\ell)$ for some semi-infinite positive braid $\B$. Hogancamp in \cite{Hog15} gives an explicit calculation of $\H(\tilde P)$ and proves that this is the stable HOMFLY-PT homology of the $(n,m)$ torus knot as $m \to \infty$.

\begin{theorem}[Hogancamp \cite{Hog15}]\label{thm-HHHofP}
There is an isomorphism of triply graded algebras $$\H(\tilde P) \simeq \mathcal{A}_n:=\QQ[u_1,\ldots,u_n,\xi_1,\ldots,\xi_n]$$ where the $u_k$ are even indeterminates of tridegree $t^{2k-2}q^{-2k}$ and the $\xi_k$ are odd indeterminates of tridegree $t^{2k-2}q^{4-2k}a$. This homology is isomorphic to a limit of HOMFLY-PT homologies of the torus links $(n,k)$ as $k\to \infty$.
\end{theorem}

\begin{remark}\label{rmk-gradings}
The gradings for the generators of $\mathcal{A}_n$, $u_k$ and $\xi_k$, are slightly different in this text than in \cite{Hog15}. However, after applying a duality functor $\mathcal{D}: \mathcal{K}^-(\mathbb{S}\text{Bim}_n) \to \mathcal{K}^+(\mathbb{S}\text{Bim}_n)$ to the projector and inverting the $q$-grading we get the form of the algebra $\mathcal{A}_n$ in this text.
\end{remark}

Since in the case of HOMFLY-PT homology, we are working with $1$-colored braids, we can sharpen the bound from the proof Theorem \ref{thm-infBraid} (Equation \ref{eq-bigHomBound}).

\begin{lemma}\label{lem-homBoundHHH}
Let $\B$ be a semi-infinite complete positive braid word on $n$ strands and $\B_\ell$ denote its partial braids. Suppose $\B_\ell$ contains $y(\ell)$ diagonals and let $F_\ell: \F(\B_\ell) \to \F(\T^{y(\ell)})$ be the canonical projection map. Then $|F_\ell|_h \geq y(\ell)$.	
\end{lemma}

\begin{proof}
In the $1$-colored case, the ladder terms in the mapping cone of $F_\ell$ all include edges labeled by zero, which can be deleted resulting in barbells (see Remark \ref{rmk-0edges}).  This in turn means that there are no rungs causing the extra complications indicated in Remark \ref{rmk-newVsOld}, and so the split and merge points of the barbell can be pulled through fractional twists in precisely the same way that the turnbacks are pulled through in the proof of Theorem 1.1 of \cite{IW16}.  The ensuing homological degree shifts lead to the same bound as in that paper, that is, the number of diagonals $y(\ell)$.
\end{proof}

\begin{corollary} \label{cor-twistEst}
	The projection map $\pi:\F(\T^{y+1}) \to \F(\T^{y})$ has homological order $|\pi|_h \geq y$. 
\end{corollary}

Lemma \ref{lem-homBoundHHH} shows that in the $1$-colored case  $|F_\ell|_h$ approaches $\infty$ at a much faster rate than the conservative estimate in the proof of Theorem \ref{thm-infBraid} which is needed for the $m$-colored case. We now describe how we can use Lemma \ref{lem-homBoundHHH} to give estimates on HOMFLY-PT homology.

\begin{theorem}\label{thm-linkEst}
Let $\b$ be a finite positive braid on $n$ strands with $\ell$ crossings. Then there exists a map $\bar{F}_\ell: \H(\b) \to \H(\T_n^{y(\ell)})$ which is an isomorphism for all homological degrees less than $y(\ell)$. 	
\end{theorem}
\begin{proof}
	Consider $\b$ as a partial braid $\b = \B_\ell$ of some positive semi-infinite braid $\B$. By Lemma \ref{lem-homBoundHHH} we know that the projection map $\F_\ell : \F(\b) \to \F(\T^{y(\ell)})$ has homological order at least $y(\ell)$. Since both $\HH$ and $H$ are functors, the result follows.
\end{proof}

Let $\mathcal{A}_{n,<y(\ell)}$ denote the subspace of $\mathcal{A}_n$ spanned by monomials $u_1^{a_1}...u_n^{a_n}\xi_1^{b_1}...\xi_n^{b_n}$ where $a_i \in \NN, b_i \in \{0,1\}$ and $\sum_{i=1}^n (2i-2)(a_i+b_i) < y(\ell)$. Combining Theorem \ref{thm-linkEst} and Theorem \ref{thm-HHHofP} we get the following corollary immediately.

\begin{corollary}\label{cor-linkEst}
	Let $\b$ be a finite positive braid on $n$ strands with $\ell$ crossings and let $\H_{<y(\ell)}(\b)$ denote the subspace of $\H(\b)$ of all terms with homological degree less than $y(\ell)$. Then 
	$$\H_{<y(\ell)}(\b) \cong \mathcal{A}_{n,<y(\ell)}.$$
\end{corollary}

\begin{proof}
Corollary \ref{cor-twistEst} and Theorem \ref{thm-linkEst} imply that the map $\H(\pi): \H(\tilde{P}) \to \H(T^{y(\ell)})$ is an isomorphism for all homological degrees less than $y(\ell)$. Composing this isomorphism with the inverse of the isomorphism (for homological degree less that $y(\ell)$) given by Theorem \ref{thm-linkEst} gives the result.
\end{proof}

\subsection{Colored HOMFLY-PT homology} \label{subsec-clrHHH}

We now briefly discuss the case of colored HOMFLY-PT homology. Colored HOMFLY-PT homology was first defined using certain bimodules by Mackaay, Stosic and Vaz \cite{MSV11}. Though they only proved invariance in the case of $1$ and $2$-colored components, they defined their complexes of bimodules in the case of arbitrary colors. Later Webster and Williamson in \cite{WW09} proved the invariance for all colors using geometric methods. The bimodules involved in the construction of colored HOMFLY-PT homology are known as singular Soergel bimodules. For the sake of brevity we will not give a rigorous definition of the entire construction but define it via analogy with the $\sln$ case.

We first associate bimodules to the basic webs from Figure \ref{fig-basicWebs}. Define $R_\ell = \QQ[x_1,...,x_\ell]^{S_\ell}$ to be the ring of symmetric polynomials in $\ell$ variables with $\deg(x_i)=-2$ as before. Also define $R_{i_1,...,i_n} = R_{i_1} \otimes_\QQ \cdots \otimes_\QQ R_{i_n}$. We associate to the basic webs in Figure \ref{fig-basicWebs} the following bimodules: 

\begin{equation}
\begin{tikzpicture}
	\draw[thick,->] (2,0)--(2,1); 
	\node at (2,-.25) {$i$};
	\node at (2,1.25) {$i$};
	\node at (3,.5) {$= R_{i}.$};
\end{tikzpicture}
\end{equation}

\begin{equation}
\begin{tikzpicture}
	\draw[thick,->-] (2,0)--(2,.5); \node at (2,-.25) {$i+j$};
	\draw[thick,->-] (2,.5)--(1.5,1); \node at (1.5,1.25) {$i$};
	\draw[thick,->-] (2,.5)--(2.5,1);\node at (2.5,1.25) {$j$};
	\node at (4,.5) {$= R_{i,j}\otimes_{R_{i,j}}R_{i+j} .$};
\end{tikzpicture}
\end{equation}

\begin{equation}
\begin{tikzpicture}
	\draw[thick,->-] (4,0)--(4.5,.5); \node at (4,-.25) {$i$};
	\draw[thick,->-] (5,0)--(4.5,.5); \node at (5,-.25) {$j$};
	\draw[thick,->-] (4.5,.5)--(4.5,1); \node at (4.5,1.25) {$i+j$};
	\node at (6.5,.5) {$= R_{i+j}\otimes_{R_{i,j}}R_{i,j} .$};
\end{tikzpicture}	
\end{equation}

In general, we can associate to every ``$\slinf$-web'' a singular Soergel bimodule. That is, we associate certain bimodules to webs with possibly arbitrarily high labels via the following rules. The monoidal product on objects of $\foamC$ corresponds to the $\otimes_\QQ$ functor on bimodules and composition of $\slinf$-webs corresponds to the tensor product over $R_{i_1,...,i_n}$ where $i_1,...,i_n$ are the intermediate labels. The relations in Figure \ref{fig-webRel} hold for $\slinf$-webs as well (see \cite{MSV11, Will11} for the exact relations).

We can now define $m$-colored HOMFLY-PT homology in an analogous manner to the definition of $m$-colored Khovanov-Rozansky homology.

Let $\B$ be a $m$-colored (finite) braid, then we define $\F(\b)$ in the same manner as $\C{\b}$, but with singular Soergel bimodules/$\slinf$-webs instead. In particular, we define 

\begin{equation}\label{eq-uniclrPosCplxHHH}
	\begin{tikzcd}
	\F\left(\posxT\right) = \arcsT \arrow{r}{d_0} & tq \ladder{\,}{\,}{\,}{\,}{1} \arrow{r}{d_1} & \cdots \arrow{r}{d_{m-2}}& t^{m-1}q^{m-1} \ladder{\,}{\,}{\,}{\,}{m-1} \arrow{r}{d_{m-1}} &t^mq^m\ladder{\,}{\,}{\,}{\,}{m}.
	\end{tikzcd}
\end{equation}
\begin{equation}\label{eq-uniclrNegCplxHHH}
	\begin{tikzcd}
	\F\left(\negxT\right) = \ladder{\,}{\,}{\,}{\,}{m}\arrow{r}{d_m'} & tq \ladder{\,}{\,}{\,}{\,}{m-1} \arrow{r}{d_{m-1}'} & \cdots \arrow{r}{d_{1}'}& t^{m-1}q^{m-1} \ladder{\,}{\,}{\,}{\,}{1} \arrow{r}{d_{0}'} &t^mq^m\arcsT.
	\end{tikzcd}
\end{equation}

The maps $d_i$ and $d_i'$ are represented by foams in the same manner as they were for colored Khovanov-Rozansky homology. However, we will omit their definition as bimodule morphisms for the sake of brevity and refer the reader to \cite{MSV11,WW09} for those details.

Our key lemmas from \S\ref{subsec-forkST} on the homological shifts with respect to fork sliding and twisting hold as well for $m$-colored HOMFLY-PT homology.

\begin{lemma}[Fork Sliding] \label{lem-forkSlideHHH}
	Let $T_1, T_1',T_2,T_2'$ be the colored tangles in Figure \ref{fig-forkSlide}.  Then 
	\begin{equation}\label{eq-forkSlideHHH}
		\begin{aligned}
			\F(T_1) &\simeq (tq)^{\min(i,k)+\min(j,k) - \min(i+j,k)}\F(T_1')\\
			\F(T_2) &\simeq \F(T_2')
		\end{aligned}
	\end{equation}
\end{lemma}

\begin{lemma}[Fork Twisting] \label{lem-forkTwistHHH}
	Let $T_3, T_3',T_4,T_4'$ be the colored tangles in Figure \ref{fig-forkTwist}.  Then 
	\begin{equation}\label{eq-forkTwistHHH}
		\begin{aligned}
			\F(T_3) &\simeq t^{\min(i,j)}q^{ij+\min(i,j)}\F(T_3'),\\
		\F(T_4) &\simeq q^{-ij}\F(T_4').
		\end{aligned}
	\end{equation}
\end{lemma}

These lemmas follow from the invariance proofs in \cite{WW09}. Combining these lemmas together, we can restate Proposition \ref{prop-isotopyShifts} for colored HOMFLY-PT homology. 

\begin{proposition}\label{prop-isotopyShiftsHHH}
Let $T_1$ and $T_2$ be any two colored positive tangles such that we can produce $T_2$ from $T_1$ by a sequence of braid-like Reidemeister moves and fork/ladder twists and slides.  Then $\F(T_1)\simeq t^\alpha q^\beta \F(T_2)$ where
\begin{equation}\label{eq-isotopyShiftsHHH}
\alpha = \sum_{\chi\in T_1} \min(\chi) - \sum_{\delta\in T_2} \min(\delta).
\end{equation}
Here each sum is taken over all crossings in each diagram.
\end{proposition}

The proof of Proposition \ref{prop-isotopyShiftsHHH} is the same as the proof of Proposition \ref{prop-isotopyShifts}.  We now present the theorem of Cautis concerning the homology of the colored infinite twist and its relation to highest weight projectors.

\begin{theorem}[Cautis \cite{Cau16}]\label{thm-infFTHHH}
Let $n \geq N$ and let $V=\CC^N_q$ denote the standard $N$-dimensional representation of $\qsln$. Then $\F(\FT_{(m)}^\infty)$ categorifies the highest weight projector $$p : (\Lambda^mV)^{\otimes n} \to (\Lambda^mV)^{\otimes n}$$ factoring through the highest weight subrepresentation of $(\Lambda^mV)^{\otimes n}$. 
\end{theorem}

\begin{theorem}\label{thm-infBraidHHH}
	Let $\BB$ be a complete semi-infinite positive braid on $n$ strands represented by the word $\B$, and define $\F(\B_{(m)}) = \lim_{\ell\to\infty} \F(\B_{\ell,(m)})$ if it exists. Then $\F(\B_{(m)})$ exists and is chain homotopy equivalent to $\F(\FT_{(m)}^\infty)$ for all $m \in \NN$, and thus we have a well-defined 
	$$ \F(\BB_{(m)}) \simeq \F(\FT_{(m)}^\infty). $$
\end{theorem}

The proof of Theorem \ref{thm-infBraidHHH} follows from the same argument as Theorem \ref{thm-infBraid}, and so we omit it here.  The corollaries in \S \ref{subsec-mainResult} hold as well for the HOMFLY-PT case, and we leave these translations to the reader.

\bibliographystyle{amsalpha}
\bibliography{bib}

\end{document}